\documentclass[12pt, a4paper]{article}

\PassOptionsToPackage{noend}{algpseudocode}
\PassOptionsToPackage{nameinlink,capitalize}{cleveref}
\PassOptionsToPackage{breakspaces}{clevethm}
\PassOptionsToPackage{inline}{enumitem}
\PassOptionsToPackage{bookmarksdepth=4}{hyperref}
\PassOptionsToPackage{dvipsnames}{xcolor}
\usepackage{tcolorbox}

\usepackage[
	alglinenumber,
	thmreset=section,
	eqreset=section,
	arxiv,
]{myPreamble}
\pgfplotsset{compat=1.16}
\usepackage{url}
\Alignfalse
\newcommand{\NL}[1][]{\ifstrempty{#1}{,}{#1}\ }

\newcommand{\keywords}[1]{\par\noindent{\def\and{\unskip,\ }{\bf Keywords. }#1}\par}
\newcommand{\subclass}[1]{\par\noindent{\def\and{\unskip,\ }{\bf AMS subject classifications. }#1}\par}
\let\OLDand\and
\def\and{\texorpdfstring{\OLDand}{, }}%

\Crefname{figure}{Figure}{Figures}

\theoremstyle{definition}
\newtheorem{remark}{Re\-mark}[section]
\setlistlabel[proofitemize]{\textbullet}

\crefname{ALG@line}{step}{steps}

\newtcolorbox{mybox}[1][]{%
	left=0pt,
	right=0pt,
	top=0pt,
	bottom=0pt,
	colback=gray!12,
	colframe=gray!12,
	width=\dimexpr\textwidth\relax,
	enlarge left by=0mm,
	boxsep=5pt,
	arc=5pt,outer arc=5pt,
	#1
}

\newcommand{\toattentive}[1]{\overset{#1}{\to}}
\DeclareMathOperator{\suchthat}{\big|} 
\newcommand{\LL}{\mathcal{L}}
\newcommand{\LLslack}{\mathcal{L}^{\text{S}}}
\newcommand{\panoc}[1]{%
	\texorpdfstring{%
		\ifstrempty{#1}{%
			PANOC%
		}{%
			\if+#1%
				PANOC$^{#1}$%
			\else
				PANOC$^{#1}$%
			\fi
		}%
	}{PANOC#1}%
}
\newcommand{\XX}{\R^n}
\newcommand{\YY}{\R^m}
\newcommand{\ZZ}{\R^p}
\newcommand{\Ybounded}{Y}
\newcommand{\normalcone}{\mathcal{N}}
\newcommand{\limnormalcone}{\normalcone^{\textup{lim}}}
\newcommand{\Deo}{{D_{\text{EO}}}}

\newcommand{\emailLink}[1]{\textsc{email} \href{mailto:#1}{#1}}
\newcommand{\orcidLink}[1]{\textsc{orcid} \href{https://orcid.org/#1}{#1}}
\newcommand{\amsmscLink}[1]{\href{http://www.ams.org/mathscinet/msc/msc2020.html?t=#1}{#1}}

\newcommand{\alps}{{ALS}}

\newcommand{\BazingaFullLink}{\href{https://github.com/aldma/Bazinga.jl}{https://github.com/aldma/Bazinga.jl}}
\newcommand{\ProximalAlgorithmsLink}{\href{https://github.com/JuliaFirstOrder/ProximalAlgorithms.jl}{ProximalAlgorithms.jl}}

\usetikzlibrary{patterns}%

\newcommand{\TheAuthorADM}{Alberto De~Marchi}
\newcommand{\TheAuthorXJ}{Xiaoxi Jia}
\newcommand{\TheAuthorCK}{Christian Kanzow}
\newcommand{\TheAuthorPM}{Patrick Mehlitz}
\newcommand{\TheTitle}{Constrained Composite Optimization and Augmented Lagrangian Methods}

\newcommand{\TheKeywords}{%
	Augmented Lagrangian methods\and
	Composite nonconvex optimization\and
	Nonlinear optimization\and
	Nonsmooth optimization
}
\newcommand{\TheAMSsubj}{%
	\amsmscLink{49J53}\and
	\amsmscLink{65K05}\and
	\amsmscLink{90C30}
}
\newcommand{\TheOrcidADM}{0000-0002-3545-6898}
\newcommand{\TheOrcidXJ}{0000-0002-7134-2169}
\newcommand{\TheOrcidCK}{0000-0003-2897-2509}
\newcommand{\TheOrcidPM}{0000-0002-9355-850X}
\newcommand{\TheAffiliationADM}{%
	Universit{\"a}t der Bun\-des\-wehr M{\"u}n\-chen\NL%
	Department of Aerospace Engineering\NL%
	Institute of Applied Mathematics and Scientific Computing\NL%
	85577 Neubiberg/Munich, Germany\NL[.]%
}
\newcommand{\TheAffiliationPM}{%
	Brandenburg University of Technology Cottbus-Senftenberg\NL%
	Institute of Mathematics\NL%
	03046 Cottbus, Germany\NL[.]%
}
\newcommand{\TheAffiliationCK}{%
	University of W{\"u}rzburg\NL%
	Institute of Mathematics\NL%
	97074 Würzburg, Germany\NL[.]%
}
\newcommand{\TheAffiliationXJ}{%
	University of W{\"u}rzburg\NL%
	Institute of Mathematics\NL%
	97074 Würzburg, Germany\NL[.]%
}
\newcommand{\TheAcknowledgements}{%
	Alberto De~Marchi is grateful to Andreas Themelis 
	(Kyushu University), for sharing his insight and rigour,
	and to Matthias Gerdts 
	(Universit{\"a}t der Bundeswehr M{\"u}nchen), for the support and guidance.
	The authors wish to thank two anonymous referees for their detailed comments and constructive suggestions,
	which significantly shaped and improved the quality of this work.
}
\newcommand{\TheFunding}{%
	Xiaoxi Jia and Christian Kanzow acknowledge support by the German Research Foundation (DFG) within the priority program
	\emph{Non-smooth and Complementarity-based Distributed Parameter Systems: 
		Simulation and Hierarchical Optimization} (SPP 1962) under grant numbers KA 1296/24-2.
}

\author{\TheAuthorADM\thanks{\TheAffiliationADM\newline\emailLink{alberto.demarchi@unibw.de}\orcidLink{\TheOrcidADM}}%
	\and\TheAuthorXJ\thanks{\TheAffiliationXJ\newline\orcidLink{\TheOrcidXJ}}%
	\and\TheAuthorCK\thanks{\TheAffiliationCK\newline\orcidLink{\TheOrcidCK}}%
	\and\TheAuthorPM\thanks{\TheAffiliationPM\orcidLink{\TheOrcidPM}}}
\date{}
\title{\TheTitle}
\begin{document}

\maketitle

\begin{abstract}
	\noindent
	We investigate finite-dimensional constrained structured optimization problems,
featuring composite objective functions and set-membership constraints.
Offering an expressive yet simple language, this problem class provides a modeling framework for a variety of applications.
We study stationarity and regularity concepts, and propose
a flexible augmented Lagrangian scheme.
We provide a theoretical characterization of the algorithm and its asymptotic properties, deriving convergence results for fully nonconvex problems.
It is demonstrated how the inner subproblems can be solved by off-the-shelf proximal methods, notwithstanding the possibility to adopt any solvers, insofar as they return approximate stationary points.
Finally, we describe our matrix-free implementation of the proposed algorithm and test it numerically.
Illustrative examples show the versatility of constrained composite programs as a modeling tool and expose difficulties arising in this vast problem class.
\end{abstract}
\subclass{\TheAMSsubj}

\clearpage

\section{Introduction}\label{sec:introduction}
	In this paper we investigate and develop numerical methods for constrained 
composite programs, namely finite-dimensional optimization problems 
of the form
\begin{equation}
    \minimize_{x}
    {}\quad{}
    q(x) \coloneqq f(x) + g(x)
    {}\qquad{}
    \stt
    {}\quad{}
    c(x) \in D ,
    \tag{P}\label{eq:P}
\end{equation}
where $x$ is the decision variable, $f$ and $c$ are smooth functions, $g$ is proper and lower semicontinuous, and $D$ is a nonempty closed set.
We call \eqref{eq:P} a constrained composite optimization problem because it contains set-membership constraints and a composite objective function $q \coloneqq f + g$.
Notice that the problem data, namely $f$, $g$, $c$ and $D$, can be nonconvex, the nonsmooth cost term $g$ can be discontinuous and the constraint set $D$ can be disconnected.
Thanks to their rich structure and flexibility, constrained composite problems are of interest for modeling in a variety of applications, ranging from optimal and model predictive control \cite{demarchi2020constrained,sopasakis2020open} to signal processing \cite{combettes2011proximal}, low-rank and sparse approximation, compressed sensing, cardinality-constrained optimization \cite{beck2018optimization} and disjunctive programming \cite{balas2018disjunctive}, such as problems with complementarity, vanishing and switching constraints \cite{jia2021augmented,mehlitz2019comparison}.

Augmented Lagrangian methods have recently attracted revived and grown interest.
Tracing back to the classical work of Hestenes \cite{hestenes1969multiplier} and Powell \cite{powell1969method}, 
the augmented Lagrangian framework can tackle large-scale constrained problems.
Recent accounts on this topic can be found in \cite{bertsekas1996constrained,birgin2014practical,conn1991globally}, among others.
Our approach is inspired by the fact that 
``augmented Lagrangian ideas are independent of the degree of smoothness of the functions that define the problem'' \cite[\S 4.1]{birgin2014practical} 
and lead to a sequence of unconstrained or simply constrained subproblems.
Moreover, this framework can handle nonconvex constraints, is often superior to pure penalty methods, 
enjoys good warm-starting capabilities and allows to avoid ill-conditioning due to a pure penalty approach 
as well as to deal with constraints without softening them; cf.\ \cite{sopasakis2020open,stella2017simple}.
In the context of constrained composite programming, 
the augmented Lagrangian subproblems associated with \eqref{eq:P} may, again, be
of composite type but possess, if at all, comparatively simple constraints. 
Exemplary, these subproblems can be solved with the aid of proximal methods, 
inaugurated by Moreau \cite{moreau1965proximite}, which can handle nonsmooth,
nonconvex and extended real-valued cost functions; cf.\ 
\cite{combettes2011proximal,KanzowMehlitz2022,parikh2014proximal,themelis2018proximal}
for recent contributions.

The close relationship between augmented Lagrangian and proximal methods is well known and traces back to Rockafellar \cite{rockafellar1976augmented}.
These approaches have been combined in \cite{dhingra2019proximal} to deal with unconstrained, composite 
optimization problems whose nonsmooth term is convex and possibly composed with a linear operator.
Following this strategy, the proximal augmented Lagrangian method has been considered for constrained composite programs 
in \cite[Ch.~1]{demarchi2021dissertation}, however lacking of sound theoretical support and convergence analysis.
A first step for resolving these shortcomings is constituted by proximal gradient methods 
that can cope with \emph{local} Lipschitz continuity of the smooth cost gradient,
only recently investigated in the Euclidean setting,
see \cite{demarchi2022proximal,KanzowMehlitz2022}.
By relying on an adaptive stepsize selection rule for the proximal gradient oracle,
these algorithms can be adopted as inner solver for augmented Lagrangian subproblems arising from general nonlinear constraints.

Another issue originates from the following observation.
One can reformulate the original problem, by introducing slack variables, in order to have a set-membership constraint with a convex right-hand side; 
consider this problem equipped with slack variables and the associated augmented Lagrangian function.
The proximal augmented Lagrangian function characterizes the latter one on the manifold 
corresponding to the explicit minimization over the slack variables \cite{dhingra2019proximal,rockafellar1976augmented}.
This procedure is employed to eliminate the slack variables and, in the convex setting, to obtain a continuously differentiable function.
Although the same ideas apply to \eqref{eq:P}, the resulting proximal augmented Lagrangian does not exhibit this favorable property in the fully nonconvex setting.
In particular, this lack of regularity is due to the set-valued projection onto the constraint set $D$.

The contribution of this work touches several aspects.
We investigate the abstract class of constrained composite optimization problems 
in the fully nonconvex setting and discuss relevant stationarity concepts.
Then, we present an algorithm for the numerical solution of these problems and, considering a classical (safeguarded) augmented Lagrangian scheme,
we provide a comprehensive yet compact global convergence analysis.
Patterning this methodology, analogous algorithms and theoretical results can be derived based on other augmented Lagrangian schemes.
Further, we demonstrate that there is no need for special choices of possibly set-valued projections and proximal mappings since we rely on the
aforementioned reformulation of \eqref{eq:P} with slack variables and keep them within our algorithmic framework. 
It is carved out that, apart from the higher number of decision variables, this reformulation is nonhazardous.
We show that it is possible to adopt off-the-shelf, yet adaptive, proximal gradient methods for solving the augmented Lagrangian subproblems.
Finally, some numerical experiments visualize computational features of our algorithmic approach.

The following blanket assumptions are considered throughout, without further mention.
Technical definitions are given in \cref{sec:preliminaries}.
\begin{mybox}
	\begin{ass}\label{ass:P}
		The following hold in \eqref{eq:P}:
		\begin{enumerate}
			\item\label{ass:f}%
				$\func{f}{\XX}{\R}$ and $\func{c}{\XX}{\YY}$ are continuously differentiable with locally Lipschitz continuous derivatives;
			\item\label{ass:g}%
				\(\func{g}{\XX}{\Rinf}\) is proper, lower semicontinuous and prox-bounded;
			\item\label{ass:set}%
				\(D\subset \YY\) is a nonempty and closed set.
		\end{enumerate}
	\end{ass}
\end{mybox}
Notice that the consequential theory remains valid whenever $\XX$ and $\YY$ are replaced by finite-dimensional Hilbert spaces $\mathbb{X}$ and $\mathbb{Y}$.
Moreover, the local Lipschitz continuity in \cref{ass:f} is actually superfluous for the 
augmented Lagrangian framework, but sufficient to solve the arising inner problems via 
proximal gradient methods \cite{demarchi2022proximal,KanzowMehlitz2022}.

By \cref{ass:f,ass:g}, the cost function $q \coloneqq f + g$ has nonempty domain, that is, $\dom q \ne \emptyset$.
Similarly, \cref{ass:set} guarantees that it is always possible to project onto the constraint set \(D\).
Nevertheless, these conditions do not imply the existence of feasible points for \eqref{eq:P};
in fact, the projection onto the set $\{ x \in \XX \,\vert\, c(x) \in D \}$ induced by the constraints $c(x) \in D$ can be as difficult as the original problem \eqref{eq:P}.
As it is the case in nonlinear programming \cite{birgin2014practical}, 
we will study the minimization properties of the augmented Lagrangian scheme with respect to some infeasibility measure.

Finally, we should mention that, for our actual implementation,
we work under the practical assumption that (only) the following computational oracles are available or simple to evaluate:
\begin{itemize}
	\item cost function value $f(x)$ and gradient $\nabla f(x)$, given $x \in \dom q$;
	\item (arbitrary) proximal point $z \in \prox_{\gamma g}(x)$ and function value $g(z)$ therein, given $x \in \XX$ and $\gamma \in (0,\gamma_g)$, $\gamma_g$ being the prox-boundedness threshold of $g$;
	\item constraint function value $c(x)$ and Jacobian-vector product $\nabla c(x)^\top v$, given $x \in \dom q$ and $v \in \YY$;
	\item (arbitrary) projected point $z \in \proj_D(v)$, given $v \in \YY$.
\end{itemize}
Relying only on these oracles, the method considered for our numerical examples is first-order and matrix-free by construction; as such, it involves only simple operations and has low memory footprint.

	\subsection{Related Work}
	Augmented Lagrangian schemes have been extensively investigated \cite{bertsekas1996constrained,birgin2014practical,conn1991globally,sopasakis2020open}, also in the infinite-dimensional setting \cite{antil2020alesqp,BoergensKanzowMehlitzWachsmuth2019,kanzow2018augmented}.

Merely lower semicontinuous cost functions have been considered in \cite{evens2021neural}. 
Inspired by \cite[Alg.\ 1]{grapiglia2020complexity} and leveraging the idea behind \cite[Ex.\ 4.12]{birgin2014practical}, 
the convergence properties of \cite[Alg.\ 1]{evens2021neural} hinge on the upper boundedness of the augmented Lagrangian 
along the iterates ensured by the initialization at a feasible point.
Although possible in some cases, in general finding a feasible starting point can be as hard as the original problem.
We deviate in this respect, seeking instead a method able to start from any $x^0 \in \XX$.
Nonetheless, if a feasible point is readily available for \eqref{eq:P}, 
one can adopt \cite[Alg.\ 1]{evens2021neural} in its original form, replacing the augmented Lagrangian function and inner solver accordingly.
In this case, and possibly assuming lower boundedness of the cost function $q$, stronger convergence guarantees can be obtained.

Programs with geometric constraints have been studied in 
\cite{BoergensKanzowMehlitzWachsmuth2019,jia2021augmented} and,
for the special case of so-called complementarity constraints, in \cite{guo2021new}.
These have a continuously differentiable cost function $f$ and set-membership constraints of the form $c(x) \in C$, $x \in D$, 
with $D$ as in \cref{ass:set} and $C$ nonempty, closed and \emph{convex}.
As already mentioned, similar structure can be obtained from \eqref{eq:P} by introducing slack variables.
Moreover, as pointed out in \cite[\S 5.4]{jia2021augmented}, considering a lower semicontinuous functional $q \coloneqq f + g$ does not enlarge the problem class, since there is an equivalent, yet smooth, reformulation in terms of the epigraph of $g$.
These observations imply that constrained composite programs do not generalize the problem class considered in \cite{jia2021augmented}.
Nevertheless, the necessary reformulations come at a price: increased problem size due to slack variables and the need for projections onto the epigraph of $g$.
The augmented Lagrangian method we are about to present is designed around \eqref{eq:P} in the fully nonconvex setting. Hence, it natively handles nonsmooth cost functions, nonlinear constraints and nonconvex sets, with no need for oracles other than those mentioned above.
Analogous considerations hold for \cite{chen2017augmented}, dedicated to an augmented Lagrangian method for non-Lipschitz nonlinear programs, 
and \cite[\S 6.2]{KrugerMehlitz2021}, where the solution of the augmented Lagrangian subproblems is not discussed.

The work presented in this paper collects and builds upon some ideas put forward in \cite{demarchi2021dissertation}.
However, we consider different stationarity concepts and necessary optimality conditions, 
not based on the proximal operator as in \cite[\S 1.2]{demarchi2021dissertation}, but rather exploiting tools from variational analysis; 
see \cite{guo2018necessary,jia2021augmented,KrugerMehlitz2021,mehlitz2020asymptotic}.
Furthermore, by avoiding the marginalization approach of \cite[\S 1.4]{demarchi2021dissertation} and so maintaining the slack variables explicit, we can offer rigorous convergence guarantees for the subproblems \cite{demarchi2022proximal,KanzowMehlitz2022}, transcending the dubious justifications given in \cite[\S 1.5.4]{demarchi2021dissertation}.

\section{Notation and Fundamentals}
	In this section, we comment on notation, preliminary definitions and useful results.
	
	\subsection{Preliminaries}\label{sec:preliminaries}
	With $\R$ and $\Rinf \coloneqq \R \cup \{\infty\}$ we denote the real and extended real line,
respectively.
Furthermore, let $\R_+$ and $\R_{++}$ be the nonnegative and positive real numbers,
respectively.
We use $0$ in order to represent the scalar zero as well as the zero vector of 
appropriate dimension.
The vector in $\R^n$ with all elements equal to  $1$ is denoted by $1_n$.
The effective domain of an extended real-valued function $\func{h}{\R^n}{\Rinf}$ is denoted by $\dom h \coloneqq \{x\in\R^n \suchthat h(x) < \infty \}$.
We say that $h$ is \emph{proper} if $\dom h \neq \emptyset$
and \emph{lower semicontinuous} (lsc) if $h(\bar{x}) \leq \liminf_{x\to\bar{x}} h(x)$ for all $\bar{x} \in \R^n$.

Given a proper and lsc function $\func{h}{\XX}{\Rinf}$ and a point $\bar{x}\in\dom h$, we may avoid to assume $h$ continuous and instead appeal 
to $h$-\emph{attentive} convergence of a sequence $\{ x^k \}$:
\begin{equation}
x^k \toattentive{h} \bar{x}
{}\quad:\Leftrightarrow\quad{}
x^k \to \bar{x}
{}\quad\text{with}\quad{}
h(x^k) \to h(\bar{x}) .
\end{equation}
Following \cite[Def. 8.3]{rockafellar1998variational}, we denote by $\ffunc{\hat{\partial} h}{\XX}{\XX}$ the \emph{regular subdifferential} of $h$, where
\begin{equation}
v \in \hat{\partial} h(\bar{x})
{}\quad:\Leftrightarrow\quad{}
\liminf_{\substack{x\to\bar{x}\\x\neq\bar{x}}} \frac{h(x) - h(\bar{x}) - \langle v, x-\bar{x}\rangle}{\|x-\bar{x}\|} \geq 0 .
\end{equation}
The (limiting) \emph{subdifferential} of $h$ is $\ffunc{\partial h}{\XX}{\XX}$, where $v \in \partial h(\bar{x})$ if and only if there exist sequences $\{x^k\}$ and $\{v^k\}$ such that $x^k \toattentive{h} \bar{x}$ and $v^k \in \hat{\partial} h(x^k)$ with $v^k \to v$.
The subdifferential of $h$ at $\bar{x}$ satisfies $\partial(h+h_0)(\bar{x}) = \partial h(\bar{x}) + \nabla h_0(\bar{x})$ for any $\func{h_0}{\XX}{\Rinf}$ continuously differentiable around $\bar{x}$ \cite[Ex. 8.8]{rockafellar1998variational}.
For formal completeness, we set $\hat{\partial} h(\bar{x})\coloneqq\partial h(\bar{x})\coloneqq\emptyset$ for each $\bar{x}\notin\dom h$.
With respect to the minimization of $h$, we say that $x^\ast \in \dom h$ is \emph{stationary} if $0 \in \partial h(x^\ast)$, which constitutes a necessary condition for the optimality of $x^\ast$ \cite[Thm 10.1]{rockafellar1998variational}.
Furthermore, we say that $x^\ast \in \XX$ is \emph{$\varepsilon$-stationary} for some $\varepsilon \geq 0$ if
\begin{equation}
	\exists \eta \in \partial h(x^\ast) \colon \|\eta\| \leq \varepsilon .
	\label{eq:epsstationary}
\end{equation}

A mapping $\ffunc{S}{\XX}{\YY}$ is \emph{locally bounded} at a point $\bar{x} \in \XX$ if for some neighborhood $V$ of $\bar{x}$ the set $S(V) \subset \YY$ is bounded \cite[Def.\ 5.14]{rockafellar1998variational}; it is called locally bounded (on $\XX$) if this holds at every $\bar{x} \in \XX$.
If $S(\bar{x})$ is nonempty, we define the \emph{outer limit} of $S$ at $\bar{x}$
by means of
\begin{align*}
	\limsup\limits_{x\to\bar{x}} S(x)
	\coloneqq
	\{
		y\in\YY
		\suchthat
		\exists x^k\to\bar{x},\,
		\exists y^k\to y,\,
		y^k\in S(x^k)\,\forall k\in\N
	\}
\end{align*}
and note that this is a closed superset of $S(\bar{x})$ by definition.

Given a parameter value $\gamma>0$, the \emph{proximal} mapping $\prox_{\gamma h}$ is defined by
\begin{equation*}
	\prox_{\gamma h}(x)
	{}\coloneqq{} 
	\argmin_{z} \left\{ h(z) + \frac{1}{2\gamma}\|z-x\|^2 \right\} ,
\end{equation*}
and we say that $h$ is \emph{prox-bounded} if it is proper and $h + \|\cdot\|^2 / (2\gamma)$ is bounded below on $\XX$ for some $\gamma > 0$.
The supremum of all such $\gamma$ is the threshold $\gamma_h$ of prox-boundedness for $h$.
In particular, if $h$ is bounded below by an affine function, then $\gamma_h = \infty$.
When $h$ is lsc, for any $\gamma \in (0,\gamma_h)$ the proximal mapping $\prox_{\gamma h}$ is locally bounded, nonempty- and compact-valued \cite[Thm 1.25]{rockafellar1998variational}.

Some tools of variational analysis will be exploited in order to describe the geometry of the nonempty, closed, but not necessarily convex set $D \subset \YY$, appearing in the formulation of \eqref{eq:P}.
The \emph{projection} mapping $\proj_D$ and the \emph{distance} function $\dist_D$ are defined by
\begin{equation*}
	\proj_D(v)
	{}\coloneqq{}
	\argmin_{z \in D} \|z - v\|
	{}\quad\text{and}\quad{}
	\dist_D(v)
	{}\coloneqq{}
	\inf_{z \in D} \|z - v\| .
\end{equation*}
The former is a set-valued mapping whenever $D$ is nonconvex, whereas the latter is always single-valued.
The \emph{indicator} function of a set $D \subset \YY$ is the function $\func{\indicator_D}{\YY}{\Rinf}$ defined as $\indicator_D(v) = 0$ if $v \in D$, and $\indicator_D(v) = \infty$ otherwise.
If $D$ is nonempty and closed, then $\indicator_D$ is proper and lsc.
The proximal mapping of $\indicator_D$ is the projection $\proj_D$; thus, $\proj_D$ is locally bounded.
Given $z \in D$, the \emph{limiting normal cone} to $D$ at $z$ is the closed cone
\begin{equation*}
	\limnormalcone_D(z) \coloneqq \limsup_{v \to z} \; \cone \left( v - \proj_D(v) \right) .
\end{equation*}
For $\tilde z\notin D$, we formally set $\limnormalcone_D(\tilde z):=\emptyset$.
The limiting normal cone is robust in the following sense:
\begin{align*}
	\limnormalcone_D(z) = \limsup_{v \to z} \; \limnormalcone_D(v).
\end{align*}
Observe that, for all $v, z \in \YY$, we have the implication
\begin{equation}
	z \in \proj_D( v )
	{}\quad\Rightarrow\quad{}
	v - z \in \limnormalcone_D(z) ,
	\label{eq:projcone}
\end{equation}
and the converse implication holds, exemplary, if $D$ is convex.
For any proper and lsc function $\func{h}{\XX}{\Rinf}$ 
and a point $\bar{x}$ with $h(\bar{x})$ finite, we have
\begin{align*}
	\partial h(\bar{x})
	=
	\left\{v\in\XX \suchthat (v,-1)\in\limnormalcone_{\epi h}(\bar{x},h(\bar{x}))\right\}
\end{align*}
where $\epi h \coloneqq \{(x,\alpha)\in\XX\times\R \suchthat h(x)\leq\alpha\}$ denotes
the epigraph of $h$.

\begin{mybox}
	\begin{lem} \label{lem:squared_distance}
		Let $D\subset\YY$ be nonempty, closed and convex.
		Furthermore, let $\func{c}{\XX}{\YY}$ be continuously differentiable.
		We consider the function $\func{\vartheta}{\XX}{\R}$ given by
		$\vartheta(x)\coloneqq\tfrac{1}{2}\dist_D^2(c(x))$ for all $x\in\XX$.
		Then, $\vartheta$ is continuously differentiable,
		and for each $\bar{x}\in\XX$, we have
		\[
			\nabla\vartheta(\bar{x})
			=
			\nabla c(\bar{x})^\top\bigl(c(\bar{x})-\proj_D(c(\bar{x}))\bigr).
		\]
	\end{lem}
\end{mybox}
\begin{proof}
	We define $\func{\psi}{\YY}{\R}$ by means of $\psi(y)\coloneqq\tfrac{1}{2}\dist_D^2(y)$ for
	all $y\in\YY$ and observe that $\vartheta=\psi\circ c$.
	Since $D$ is assumed to be convex, $\psi$ is continuously differentiable with gradient
	$
	\nabla\psi(\bar{y})=\bar{y}-\proj_D(\bar{y})
	$
	for each $\bar{y}\in\YY$,
	see \cite[Cor.\ 12.30]{BauschkeCombettes2011}, and the
	statements of the lemma follow trivially from the standard chain rule.
\end{proof}

	\subsection{Stationarity Concepts and Qualification Conditions}\label{sec:concepts}
	We now define some basic concepts and discuss stationarity conditions for \eqref{eq:P}.
As the cost function $q \coloneqq f+g$ is possibly extended real-valued, 
feasibility of a point must account for its domain.
\begin{mybox}
	\begin{defin}[Feasibility]\label{def:feasible}
		A point $x^\ast \in \XX$ is called \emph{feasible} for \eqref{eq:P} 
		if $x^\ast \in \dom q$ and $c(x^\ast) \in D$.
	\end{defin}
\end{mybox}
Working under the assumption that the constraint set $D$ is nonconvex, 
a plausible stationarity concept for addressing \eqref{eq:P} is that of 
Mordukhovich-stationarity, which exploits limiting normals to $D$; cf.\ 
\cite[\S 3]{mehlitz2020asymptotic} and \cite[Thm~5.48]{mordukhovich2006variational}.
\begin{mybox}
	\begin{defin}[M-stationarity]\label{def:Mstationary}
		Let $x^\ast \in \XX$ be a feasible point for \eqref{eq:P}.
		Then, $x^\ast$ is called a \emph{Mordukhovich-stationary} point of \eqref{eq:P} 
		if there exists a multiplier $y^\ast \in \YY$ such that
		\begin{subequations}
			\begin{align}
				- \nabla c(x^\ast)^\top y^\ast
				{}\in{}&
				\partial q(x^\ast) ,
				\label{eq:Mstationary:x} \\
				y^\ast
				{}\in{}&
				\limnormalcone_D( c(x^\ast) ) .
				\label{eq:Mstationary:y}
			\end{align}
		\end{subequations}
	\end{defin}
\end{mybox}
Notice that these conditions implicitly require the feasibility of $x^\ast$, 
for otherwise the subdifferential and limiting normal cone would be empty.
Note that this definition coincides with the usual KKT conditions of \eqref{eq:P} 
if $g$ is smooth and $D$ is a convex set.

Subsequently, we study an asymptotic counterpart of this definition. 
In case where $q$ is locally Lipschitz continuous, one could apply the notions from
\cite[\S 2.2]{jia2021augmented} and \cite[\S 5.1]{mehlitz2020asymptotic} 
for that purpose. 
However, since $g$ is assumed to be merely lsc,
we need to adjust these concepts at least slightly.
\begin{mybox}
	\begin{defin}[AM-stationarity]\label{def:AMstationary}
		Let $x^\ast \in \XX$ be a feasible point for \eqref{eq:P}.
		Then, $x^\ast$ is called an \emph{asymptotically M-stationary} point of \eqref{eq:P} 
		if there exist sequences $\{x^k\}, \{\eta^k\} \subset \XX$ and 
		$\{y^k\}, \{\zeta^k\} \subset \YY$ such that 
		$x^k \toattentive{q} x^\ast$, $\eta^k \to 0$, $\zeta^k \to 0$ and
		\begin{subequations}
			\begin{align}
				-\nabla c(x^k)^\top y^k + \eta^k
				{}\in{}&
				\partial q(x^k) ,
				\label{eq:AMstationary:x} \\
				y^k 
				{}\in{}&
				\limnormalcone_D( c(x^k) + \zeta^k )
				\label{eq:AMstationary:y}
			\end{align}
			\label{eq:AMstationary}
		\end{subequations}
		for all $k \in \N$.
	\end{defin}
\end{mybox}

The definition of an AM-stationary point is similar to the notion 
of an asymptotic KKT point \cite{birgin2014practical}, 
as well as the meaning of the iterates $x^k$ and the Lagrange multipliers $y^k$.
Notice that \cref{def:AMstationary} does not require the sequence $\{ y^k \}$ to converge. 
The vector $\eta^k$ measures the dual infeasibility,
namely the inexactness in the stationarity condition \eqref{eq:AMstationary:x} 
at $x^k$ and $y^k$.
The vector $\zeta^k$ is introduced to account for the fact that the condition 
$c(x^k) \in D$ can be violated along the iterates, though it (hopefully) holds asymptotically.
As the corresponding (limiting) normal cone $\limnormalcone_D(c(x^k))$ would be empty 
in this case, it would not be possible to satisfy the inclusion 
$y^k \in \limnormalcone_D(c(x^k))$.
The sequence $\{\zeta^k\}$ remedies this issue and gives 
a measure of primal infeasibility, as we will attest.
Finally, the convergence $x^k\toattentive{q} x^*$, which is not restrictive in situations where
$g$ is continuous (relative to its domain), will be important later on when taking the limit in
\eqref{eq:AMstationary:x} since we aim to recover the limiting subdifferential
of the objective function as stated in \eqref{def:Mstationary}.
Let us note that a slightly different notion of asymptotic stationarity has been introduced
for rather general optimization problems in Banach spaces 
in \cite[Def.~6.4, Rem.~6.5]{KrugerMehlitz2021}. 
Therein, different primal sequences are used for the objective function and the
constraints.

A local minimizer for \eqref{eq:P} is M-stationary only under validity 
of a suitable qualification condition, which, by non-Lipschitzness of
$g$, will depend on the latter function as well,
see \cite{guo2018necessary} for a discussion.
However, we can show that each local minimizer of \eqref{eq:P} is always AM-stationary.
Related results can be found in
\cite[Thm~6.2]{KrugerMehlitz2021} and \cite[\S 5.1]{mehlitz2020asymptotic}.
\begin{mybox}
	\begin{prop}\label{prop:minimizers_AM_stationary}
		Let $x^*\in\XX$ be a local minimizer for \eqref{eq:P}.
		Then, $x^*$ is an AM-stationary point for \eqref{eq:P}.
	\end{prop}
\end{mybox}
\begin{proof}
	By local optimality of $x^\ast$ for \eqref{eq:P}, we find some
	$\varepsilon>0$ such that $q(x)\geq q(x^\ast)$ is valid for
	all $x\in \mathbb{B}_\varepsilon(x^\ast) \coloneqq \left\{x\in\XX \suchthat \|x-x^\ast\|\leq\varepsilon\right\}$
	which are feasible for \eqref{eq:P}.
	Consequently, $x^\ast$ is the uniquely determined global minimizer
	of
	\begin{equation}\label{eq:regularized_problem}
		\begin{aligned}
			&\minimize_{x}
			{}\quad{}&
			&q(x)+\frac12\|x-x^\ast\|^2&
			\\
			&\stt
			{}\quad{}&
			&c(x) \in D,
			{}\quad{}
			x\in\mathbb{B}_\varepsilon(x^\ast).&
		\end{aligned}
	\end{equation}
	Let us now consider the penalized surrogate problem
	\begin{equation}\label{eq:penalized_surrogate}\tag{P$(k)$}
		\begin{aligned}
			&\minimize_{x,s}
			{}\quad{}&
			&q(x)+\frac{k}2\|c(x)-s\|^2+\frac12\|x-x^\ast\|^2&
			\\
			&\stt
			{}\quad{}&
			&x\in\mathbb{B}_\varepsilon(x^\ast),
			{}\quad{}
			s\in D\cap\mathbb{B}_1(c(x^\ast))&
		\end{aligned}
	\end{equation}
	where $k\in\N$ is arbitrary.
	Noting that the objective function of this optimization problem
	is lsc while its feasible set is nonempty and
	compact, it possesses a global minimizer $(x^k,s^k)\in\XX\times\YY$
	for each $k\in\N$. Without loss of generality, we assume
	$x^k\to\tilde{x}$ and $s^k\to\tilde{s}$ for some
	$\tilde{x}\in\mathbb{B}_\varepsilon(x^\ast)$ and $\tilde{s}\in D\cap\mathbb{B}_1(c(x^\ast))$.
	
	We claim that $\tilde{x}=x^\ast$ and $\tilde{s}=c(x^\ast)$.
	To this end, we note that $(x^\ast,c(x^\ast))$ is feasible to \eqref{eq:penalized_surrogate}
	which yields the estimate
	\begin{equation}\label{eq:estimate_from_penalization}
		q(x^k)+\frac{k}{2}\|c(x^k)-s^k\|^2+\frac12\|x^k-x^\ast\|^2
		\leq
		q(x^\ast)
	\end{equation}
	for each $k\in\N$.
	Using lower semicontinuity of $q$ as well as the convergences $c(x^k)\to c(\tilde{x})$
	and $s^k\to\tilde{s}$, taking the limit $k\to\infty$ in
	\eqref{eq:estimate_from_penalization} gives $c(\tilde{x})=\tilde{s}\in D$.
	Particularly, $\tilde{x}$ is feasible for \eqref{eq:regularized_problem}.
	Therefore, the local optimality of $x^\ast$ implies $q(x^\ast)\leq q(\tilde{x})$.
	Furthermore, we find
	\begin{multline*}
		q(\tilde{x})+\frac12\|\tilde{x}-x^\ast\|^2
		\leq
		\liminf\limits_{k\to\infty}
		\left(
		q(x^k)+\frac{k}{2}\|c(x^k)-s^k\|^2+\frac12\|x^k-x^\ast\|^2
		\right)
		\\\leq
		q(x^\ast)
		\leq
		q(\tilde{x}).
	\end{multline*}
	Hence, $\tilde{x}=x^\ast$, and noting that \eqref{eq:estimate_from_penalization}
	gives $q(x^k)\leq q(x^\ast)$ for each $k\in\N$,
	\begin{align*}
		q(x^\ast)
		\leq
		\liminf\limits_{k\to\infty}q(x^k)
		\leq
		\limsup\limits_{k\to\infty}q(x^k)
		\leq
		q(x^\ast),
	\end{align*}
	\ie, $x^k\toattentive{q}x^\ast$ follows.
	
	Due to $x^k\to x^\ast$ and $s^k\to c(x^\ast)$, we may assume without loss of
	generality that $\{x^k\}$ and $\{s^k\}$ are taken from the interior
	of $\mathbb{B}_\varepsilon(x^\ast)$ and $\mathbb{B}_1(c(x^\ast))$, respectively.
	Thus, for each $k\in\N$, $(x^k,s^k)$ is an unconstrained local minimizer of
	\[
	(x,s)\mapsto q(x)+\frac{k}2\|c(x)-s\|^2+\frac12\|x-x^\ast\|^2+\delta_D(s).
	\]
	Let us introduce $\func{\theta}{\XX\times\YY}{\Rinf}$
	by means of $\theta(x,s) \coloneqq g(x)+\indicator_D(s)$ for each pair $(x,s)\in\XX\times\YY$.
	Applying \cite[Prop.~1.107 and 1.114]{mordukhovich2006variational}, we
	find
	\[
	(0,0)
	\in 
	\bigl(\nabla f(x^k)+k\,\nabla c(x^k)^\top(c(x^k)-s^k)+x^k-x^\ast,k(s^k-c(x^k)\bigr)
	+
	\partial\theta(x^k,s^k)
	\]
	for each $k\in\N$.
	The decoupled structure of $\theta$ and \cite[Thm~3.36]{mordukhovich2006variational}
	yield the inclusion
	$\partial\theta(x^k,s^k)\subset\partial g(x^k)\times\limnormalcone_D(s^k)$
	for each $k\in\N$.
	Thus, setting $\eta^k \coloneqq x^\ast-x^k$, $y^k \coloneqq k(c(x^k)-s^k)$ and $\zeta^k \coloneqq s^k-c(x^k)$
	for each $k\in\N$ while observing that
	$\partial q(x^k)=\nabla f(x^k)+\partial g(x^k)$ holds, we have shown
	that $x^\ast$ is AM-stationary for \eqref{eq:P}.
\end{proof}

In order to guarantee that local minimizers for \eqref{eq:P} are not only
AM- but already M-stationary, the presence of a qualification condition
is necessary. The subsequent definition generalizes 
the constraint qualification from \cite[\S 3.2]{mehlitz2020asymptotic} to the
non-Lipschitzian setting and is closely
related to the so-called \emph{uniform qualification condition} 
introduced in \cite[Def.~6.8]{KrugerMehlitz2021}.
\begin{mybox}
	\begin{defin}[AM-regularity]\label{def:AMregularity}
		Let $x^\ast \in \XX$ be a feasible point for \eqref{eq:P}.
		Define the set-valued mapping $\ffunc{\mathcal{M}}{\XX \times \YY}{\XX}$ by
		\begin{equation*}
			\mathcal{M}(x,z) \coloneqq \partial g(x)+\nabla c(x)^\top \limnormalcone_D(c(x) - z) .
		\end{equation*}
		Then, $x^\ast$ is called \emph{asymptotically M-regular} for \eqref{eq:P} if
		\begin{equation*}
			\limsup_{\substack{x \toattentive{g} x^\ast\\z \to 0\phantom{^\ast}}} 
			\mathcal{M}(x,z) \subset \mathcal{M}(x^\ast,0) .
		\end{equation*}
	\end{defin}
\end{mybox}
Let us point the reader's attention to the fact that AM-regularity is not a
constraint qualification for \eqref{eq:P} in the narrower sense since it
depends explicitly on the objective function.
However, note that AM-regularity of some feasible point $x^\ast\in\XX$ for
\eqref{eq:P} reduces to
\begin{equation}\label{eq:AM_regularity_Lipschitzian_setting}
	\limsup_{\substack{x \to x^\ast\\z \to 0\phantom{^\ast}}} 
		\nabla c(x)^\top \limnormalcone_D(c(x) - z) 
		\subset 
		\nabla c(x^\ast)^\top \limnormalcone_D(c(x^\ast))
\end{equation}
whenever $g$ is locally Lipschitz continuous around $x^\ast$
since $x\rightrightarrows\partial g(x)$ is locally bounded 
at $x^\ast$ in this case, see \cite[Cor.~1.81]{mordukhovich2006variational}.
We also observe that \eqref{eq:AM_regularity_Lipschitzian_setting} corresponds
to the concept of AM-regularity which has been used in
\cite{jia2021augmented,mehlitz2020asymptotic} where $q$ is assumed to be
at least locally Lipschitz continuous, and this condition has been shown
to serve as a comparatively weak constraint qualification. 
Sufficient conditions for the validity of the more general qualification
condition from \cref{def:AMregularity} can be distilled in a similar way as
in \cite{KrugerMehlitz2021}. 
 
As a corollary of \cref{prop:minimizers_AM_stationary}, we find the following
result, along the lines of \cite[Prop.~6.9]{KrugerMehlitz2021}.

\begin{mybox}
	\begin{cor}\label{cor:AMregular_points}
		Let $x^\ast\in\XX$ be an AM-regular AM-stationary point for \eqref{eq:P}.
		Then, $x^\ast$ is an M-stationary point for \eqref{eq:P}.
		Particularly, each AM-regular local minimizer for \eqref{eq:P} is
		M-stationary.
	\end{cor}
\end{mybox}

Following the lines of the proofs of
\cite[Thm~3.2]{AndreaniMartinezRamosSilva2016} or
\cite[Thm~4.6]{BoergensKanzowMehlitzWachsmuth2019},
it is even possible to show that whenever, for each continuously differentiable 
function $f$, AM-stationarity of a feasible point $x^\ast\in\XX$ of
\eqref{eq:P} already implies M-stationarity of $x^\ast$, then $x^\ast$ must be
AM-regular. Relying on the terminology coined in 
\cite{AndreaniMartinezRamosSilva2016}, this means that AM-regularity is the
weakest \emph{strict} qualification condition associated with AM-stationarity.

\section{Augmented Lagrangian Method}\label{sec:ALM}
	Constrained minimization problems such as \eqref{eq:P} are amenable to be addressed by means of augmented Lagrangian methods.
Introducing the slack variable $s \in \YY$, \eqref{eq:P} can be rewritten as
\begin{equation}
	\minimize_{x ,\; s}
	{}\quad{}
	q(x)
	{}\qquad{}
	\stt
	{}\quad{}
	c(x) - s = 0
	,\quad
	s \in D .
	\tag{P$_{\text{S}}$}\label{eq:ALMz:P}
\end{equation}
Notice that \eqref{eq:ALMz:P} is a particular problem in the form of \eqref{eq:P}.
Moreover, if $g$ is smooth, and thus so is $q$, then \eqref{eq:ALMz:P} falls into the problem class analyzed in \cite{jia2021augmented}.
Note that $x^\ast\in\XX$ is a global (local) minimizer of \eqref{eq:P} if and only if 
$(x^\ast,c(x^\ast))$ is a global (local) minimizer of \eqref{eq:ALMz:P}.
Similarly, the M-stationary points of \eqref{eq:P} and \eqref{eq:ALMz:P} correspond to each other.
An elementary calculation additionally reveals that even the AM-stationary points of \eqref{eq:P} and \eqref{eq:ALMz:P}
can be identified with each other.
\begin{mybox}
	\begin{lem}
		A feasible point $x^\ast\in \XX$ of \eqref{eq:P} is AM-stationary for \eqref{eq:P}
		if and only if $(x^\ast,c(x^\ast))$ is AM-stationary for \eqref{eq:ALMz:P}.
	\end{lem}
\end{mybox}
\begin{proof}
	The implication $\Rightarrow$ is obvious, so let us only prove the converse one.
	If $(x^\ast,c(x^\ast))$ is AM-stationary for \eqref{eq:ALMz:P}, we find sequences
	$\{x^k\},\{\eta^k_1\}\subset\XX$ and $\{s^k\},\{y^k_1\},\{y^k_2\},\{\eta^k_2\},\{\zeta^k_1\},\{\zeta^k_2\}\subset\YY$
	such that $x^k \toattentive{q} x^\ast$, $s^k \to c(x^\ast)$, $\eta^k_i\to 0$, $\zeta^k_i\to 0$, $i=1,2$, and
	\begin{subequations}\label{eq:AM_St_with_slack}
		\begin{align}
			\label{eq:AM_St_with_slacks_x}
			-\nabla c(x^k)^\top y^k_1+\eta^k_1&\in\partial q(x^k),\\
			\label{eq:AM_St_with_slacks_s}
			y^k_1-y^k_2+\eta^k_2&=0,\\
			\label{eq:AM_St_with_slacks_normals_=}
			c(x^k)-s^k+\zeta^k_1&=0,\\
			\label{eq:AM_St_with_slacks_normals_D}
			y^k_2&\in\limnormalcone_D(s^k+\zeta^k_2)
		\end{align}
	\end{subequations}
	for all $k\in\N$, where we already used the Cartesian product rule for the limiting normal cone, cf.\ \cite[Prop.\ 1.2]{mordukhovich2006variational},
	in order to split
	\[
		(y^k_1,y^k_2)\in\limnormalcone_{\{0\}\times D}(c(x^k)-s^k+\zeta^k_1,s^k+\zeta^k_2)
	\]
	into \eqref{eq:AM_St_with_slacks_normals_=} and \eqref{eq:AM_St_with_slacks_normals_D}.
	Now, for each $k\in \N$, set $y^k:=y^k_2$, $\eta^k:=\nabla c(x^k)^\top\eta^k_2+\eta^k_1$ and $\zeta^k:=s^k-c(x^k)+\zeta^k_2$.
	Then, \eqref{eq:AMstationary:x} follows from \eqref{eq:AM_St_with_slacks_x} and \eqref{eq:AM_St_with_slacks_s}.
	Furthermore, \eqref{eq:AMstationary:y} can be distilled from \eqref{eq:AM_St_with_slacks_normals_D}.
	The convergence $\eta^k\to 0$ is clear from continuous differentiability of $c$,
	and $\zeta^k\to 0$ follows from $c(x^k)-s^k\to 0$ which is a consequence of the continuity of $c$ (or \eqref{eq:AM_St_with_slacks_normals_=}).
	\qedhere
\end{proof}

Summarizing the above observations, the way we incorporated the slack variable in \eqref{eq:ALMz:P} does not change the solution
and stationarity behavior when compared with \eqref{eq:P}. In light of \cite{BenkoMehlitz2021b}, where similar issues are
discussed in a much broader context, this is remarkable.
We use the lifted reformulation \eqref{eq:ALMz:P} as a theoretical tool to develop our approach for solving \eqref{eq:P} and investigate its properties.
For some penalty parameter $\mu>0$,
let us define the $\mu$-augmented Lagrangian function $\func{\LLslack_\mu}{\XX \times \YY \times \YY}{\Rinf}$ associated to \eqref{eq:ALMz:P} as
\begin{align}
	\LLslack_\mu(x,s,y)
	{}\coloneqq{}&
	q(x) + \indicator_D(s) + \innprod{y}{c(x) - s} + \frac{1}{2\mu} \| c(x) - s \|^2 \nonumber \\
	{}={}&
	q(x) + \indicator_D(s) + \frac{1}{2\mu} \| c(x) + \mu y - s \|^2 - \frac{\mu}{2}\|y\|^2 .
	\label{eq:Lz}
\end{align}
Observe that, by adopting the indicator $\indicator_D$, the constraint $s \in D$ is considered hard, in the sense that it must be satisfied exactly.
These simple, nonrelaxable lower-level constraints have been discussed, \eg, in \cite{andreani2008augmented,birgin2014practical,conn1991globally,jia2021augmented}.
For later use, let us compute the subdifferential of $\LLslack_\mu$ with respect to the
variables $x$ and $s$:
\begin{subequations}
	\label{eq:dLS}
	\begin{align}
		\partial_x \LLslack_\mu(x,s,y) 
		{}={}&
		\partial q(x) + \frac{1}{\mu} \nabla c(x)^\top [c(x) + \mu y - s] ,
		\label{eq:dxLS} \\
		\partial_s \LLslack_\mu(x,s,y) 
		{}={}&
		\limnormalcone_D(s) - \frac{1}{\mu} [c(x) + \mu y - s] .
		\label{eq:dzLS}
	\end{align}
\end{subequations}

The algorithm we are about to present requires, at each inner iteration, 
the (approximate) minimization of $\LLslack_\mu(\cdot,\cdot,y)$, 
given some $\mu > 0$ and $y \in \YY$, 
while in each outer iteration, $\mu$ and $y$ are updated.
This nested-loops structure naturally arises in the augmented Lagrangian framework, 
as it does more generally in nonlinear programming.

A similar method can be obtained by exploiting the structure arising from the original problem \eqref{eq:P} in order to eliminate the slack variable $s$, on the vein of the proximal augmented Lagrangian approach \cite{demarchi2021dissertation,dhingra2019proximal}.
Given some $\mu > 0$, $x \in \XX$ and $y \in \YY$, 
the explicit minimization of $\LLslack_\mu(x,\cdot,y)$ is readily obtained and yields a set-valued mapping:
\begin{equation}
	\argmin_s \LLslack_\mu(x,s,y)
	{}={}
	\proj_D \left( c(x) + \mu y \right) .
	\label{eq:z}
\end{equation}
Evaluating the augmented Lagrangian on the set corresponding to the explicit minimization over the slack variable $s$, we obtain the (single-valued) augmented Lagrangian function $\func{\LL_\mu}{\XX \times \YY}{\Rinf}$ associated to \eqref{eq:P}:
\begin{equation}
	\LL_\mu(x,y)
	{}\coloneqq{}
	\min_s \LLslack_\mu(x,s,y)
	{}={}
	q(x) + \frac{1}{2\mu} \dist_D^2 ( c(x) + \mu y ) - \frac{\mu}{2}\|y\|^2 .
	\label{eq:L}
\end{equation}
Then, one may consider replacing the minimization of $\LLslack_\mu(\cdot,\cdot,y)$ with that of $\LL_\mu(\cdot,y)$.
Following the lines of \cite[\S 4.1]{BenkoMehlitz2021b},
one can easily check that the problems $\min \LL_\mu(\cdot,y)$ and $\min \LLslack_\mu(\cdot,\cdot,y)$ are equivalent in the sense that $x^\ast$ is a local (global) minimizer of $\min \LL_\mu(\cdot,y)$ if and only if $(x^\ast,s^\ast)$, for 
each
$s^\ast \in \argmin \LLslack_\mu(x^\ast,\cdot, y)$,
is a local (global) minimizer of $\mathcal L_\mu^S(\cdot,\cdot,y)$; cf.\ \eqref{eq:z}.
However, we highlight that the term $\func{\dist_D^2}{\YY}{\R}$ is not continuously differentiable in general, as the projection onto $D$ is a set-valued mapping, thus making this approach difficult in practice.

\begin{remark}\label{rem:convex_rhs}
	Whenever $D$ is a convex set, the augmented Lagrangian function $\LL_\mu$ from
	\eqref{eq:L} is a continuously differentiable function with a locally
	Lip\-schitz continuous gradient; cf.\ \cref{lem:squared_distance}.
	Following the literature, see e.g.\ \cite{andreani2008augmented,birgin2014practical,dhingra2019proximal}, 
	one can directly augment the corresponding set-membership constraints within the
	corresponding augmented Lagrangian framework without the need of an additional 
	slack variable. 
	In practical implementations of an augmented Lagrangian framework addressing \eqref{eq:P},
	it is, thus, recommendable to treat only the difficult set-membership constraints with a
	nonconvex right-hand side with the aid of the lifting approach discussed here.
	The remaining set-membership constraints can either be augmented without slacks 
	or remain explicitly in the constraint set of the augmented Lagrangian subproblems
	if simple enough (like box constraints).
\end{remark}

The following \cref{sec:ALM:algorithm} contains a detailed statement of our algorithmic framework, whose convergence analysis is presented in \cref{sec:ALM:convergence}.
Then, suitable termination criteria are discussed in \cref{sec:ALM:termination}.
In \cref{sec:inner} we consider the numerical solution of the augmented Lagrangian subproblems.

	\subsection{Algorithm}\label{sec:ALM:algorithm}
	This section presents an augmented Lagrangian method for the solution 
of constrained composite programs of the form \eqref{eq:P}, under \cref{ass:P}.
As the augmented Lagrangian constitutes a framework, rather than a single algorithm, 
several methods have been presented in the past decades, 
expressing the foundational ideas in different flavors.
Some prominent contributions are those in 
\cite{bertsekas1996constrained,birgin2014practical,conn1991globally,grapiglia2020complexity,kanzow2018augmented,sopasakis2020open},
and for primal-dual methods \cite{gill2012primal}.
In the following, we focus on a safeguarded augmented Lagrangian scheme inspired by 
\cite[Alg.\ 4.1]{birgin2014practical} and investigate its convergence properties.
Compared to the classical augmented Lagrangian or multiplier penalty approach 
for the solution of nonlinear programs \cite{bertsekas1996constrained}, 
this variant uses a safeguarded update rule for the Lagrange multipliers 
and has stronger global convergence properties.
Although we restrict our analysis to this specific algorithm, 
analogous results can be obtained for others with minor changes.
The overall method is stated in \cref{alg:ALM} and corresponds to
the popular augmented Lagrangian solver {Algencan} from
\cite{andreani2008augmented} applied to \eqref{eq:ALMz:P}.
Let us mention, however, that the analysis in \cite{andreani2008augmented}
does neither cover composite objective functions $q \coloneqq f + g$
nor constraints of the form $c(x) \in D$ with potentially nonconvex constraint set $D$.

\begin{algorithm}[t]
	\caption{Augmented Lagrangian method for \eqref{eq:P}}%
	\label{alg:ALM}%
	\begin{algorithmic}[1]%
	\linespread{1.2}\selectfont%
	%
	%
	\Initialize Select $\mu_0 > 0$, $\theta,\kappa \in(0,1)$ and $\Ybounded \subset \YY$ nonempty bounded
	\item[For \(k=0,1,2\ldots\)]
	\State\label{state:ALM:init}\label{state:ALM:ysafe}%
		Select $\hat{y}^k \in \Ybounded$ and $\varepsilon_k \geq 0$
	\State\label{state:ALM:subproblem}
		Compute an $\varepsilon_k$-stationary point $(x^k,s^k)\in\XX\times D$ of $\LL_{\mu_k}^{\text{S}}(\cdot,\cdot,\hat{y}^k)$
	\State\label{state:ALM:y}%
		Set $y^k \gets \hat{y}^k + [ c(x^k) - s^k ] / \mu_k$
	\If{$k = 0$ or $\| c(x^k) - s^k \| \leq \theta \, \| c(x^{k-1}) - s^{k-1} \|$}\label{state:ALM:if}
		\State\label{state:ALM:ifpassed}%
			Set $\mu_{k+1} \gets \mu_k$
	\Else
		\State\label{state:ALM:iffailed}%
			Select $\mu_{k+1} \in (0, \kappa \mu_k]$
	\EndIf
\end{algorithmic}
\end{algorithm}

First of all, a primal-dual starting point is not explicitly required.
In practice, however, the subproblems at \cref{state:ALM:subproblem} 
should be solved starting from the current primal estimate $x^{k-1}$
paired with some $s^{k-1}$, preferably an element of $\proj_D(c(x^{k-1}) + \mu_k \hat{y}^k)$ as suggested by \eqref{eq:z},
thus exploiting initial guesses.
The safeguarded dual estimate $\hat{y}^k$ is drawn from a bounded set 
$\Ybounded \subset \YY$ at \cref{state:ALM:ysafe}.
Although not necessary, the choice of $\hat{y}^k$ should also depend 
on the current dual estimate $y^{k-1}$.
Moreover, the choice of $\Ybounded$ can take advantage of \emph{a priori} knowledge of $D$ 
and its structure, in order to generate better dual estimates.
For instance, if $D \subset \YY$ is compact and convex, 
we may select $\Ybounded = [-y_{\min},y_{\max}]^m$ for some $y_{\min}, y_{\max} > 0$, 
whereas if $D = \YY_+$, we may more accurately choose $\Ybounded = [-y_{\min},0]^m$; 
cf.\ \cite{jia2021augmented,sopasakis2020open}.
In practice, it is advisable to choose the safeguarded multiplier estimate $\hat{y}^k$ 
as the projection of the Lagrange multiplier $y^{k-1}$ onto $\Ybounded$, 
thus effectively adopting the classical approach as long as $y^{k-1}$ 
remains within $\Ybounded$.

The augmented Lagrangian functions and subproblems discussed above 
appear at \cref{state:ALM:subproblem}.
\cref{sec:inner} is devoted to the numerical solution of the subproblems, 
discussing several approaches.
The subproblems are usually solved only approximately, in some sense, for the sake of computational efficiency.
More precisely, the subproblem solver needs to be able to find $\varepsilon$-stationary
points of $\LLslack_\mu(\cdot,\cdot,y)$ for arbitrarily small $\varepsilon>0$,
$\mu>0$ and $y\in\Ybounded$.

\Cref{state:ALM:y} entails the classical first-order Lagrange multiplier estimate.
The update rule is designed around \eqref{eq:dxLS} and leads to the inclusion
\eqref{eq:AMstationary:x} for the primal-dual estimate $(x^k,y^k)$.
The monotonicity test at \cref{state:ALM:if} is adopted to monitor primal 
infeasibility along the iterates.
The penalty parameter is reduced at \cref{state:ALM:iffailed} in case of insufficient decrease, 
effectively implementing a simple feedback strategy to drive $\|c(x^k) - s^k\|$ to zero.

Before proceeding to the convergence analysis, 
we highlight a different interpretation of the method.
As first observed in \cite{rockafellar1976augmented}, 
the augmented Lagrangian method on the primal problem 
has an associated proximal point method on the dual problem.
Introducing the auxiliary variable $r \in \YY$, 
we rewrite the augmented Lagrangian subproblem $\min \LLslack_\mu(\cdot,\cdot,y)$ as
\begin{equation*}
	\minimize_{x ,\; s ,\; r}
	{}\quad{}
	q(x) + \indicator_D(s) + \frac{1}{2\mu}\|r - \mu y\|^2
	{}\quad{}
	\stt
	{}\quad{}
	c(x) - s + r = 0 
\end{equation*}
and then, by eliminating the slack variable $s$, as
\begin{equation*}
	\minimize_{x ,\; r}
	{}\quad{}
	q(x) + \frac{1}{2\mu}\|r - \mu y\|^2
	{}\quad{}
	\stt
	{}\quad{}
	c(x) + r \in D .
\end{equation*}
The latter reformulation amounts to a proximal dual regularization of \eqref{eq:P} 
and corresponds to a lifted representation of $\min \LL_\mu(\cdot,y)$, 
where $\LL_\mu$ is given in \eqref{eq:L},
thus showing that the approach effectively consists in solving a sequence of subproblems, 
each one being a proximally regularized version of \eqref{eq:P}.
Yielding feasible and more regular subproblems, 
this (proximal) regularization strategy has been explored 
and exploited in different contexts; some recent works are, \eg,
\cite{demarchi2022qpdo,ma2018stabilized,potschka2021sequential}.
	
	\subsection{Convergence Analysis}\label{sec:ALM:convergence}
	Throughout our convergence analysis, we assume that \cref{alg:ALM} is well-defined, 
thus requiring that each subproblem at \cref{state:ALM:subproblem} admits an approximate stationary point.
Moreover, the following statements assume the existence of some accumulation point $x^\ast$ 
or $(x^\ast,s^\ast)$
for a sequence $\{ x^k \}$ or $\{(x^k,s^k)\}$, respectively, generated by \cref{alg:ALM}.
In general, coercivity or (level) boundedness arguments should be adopted 
to verify this precondition; cf.\ \cref{lem:ALM:bounded} as well.

Due to their practical importance, we focus on affordable, or \emph{local}, 
solvers, which return merely stationary points, for the subproblems at 
\cref{state:ALM:subproblem}.
Instead, we do not present results on the case 
where the subproblems are solved to \emph{global} optimality.
The analysis would follow the classical results in \cite[Ch.~5]{birgin2014practical} 
and \cite{kanzow2018augmented}, see \cite[\S 6.2]{KrugerMehlitz2021} as well.
In summary, feasible problems would lead to feasible accumulation points 
that are global minima, in case of existence.
For infeasible problems, infeasibility would be minimized and the objective 
cost minimum for the minimal infeasibility.

Like all penalty-type methods in the nonconvex setting, 
\cref{alg:ALM} may generate accumulation points that are infeasible for \eqref{eq:P}.
Patterning standard arguments, the following result gives conditions that guarantee feasibility of limit points;
cf.\ \cite[Ex.\ 4.12]{birgin2012augmented}, \cite[Prop.\ 4.1]{jia2021augmented}.

\begin{mybox}
	\begin{prop}\label{lem:ALM:bounded}
		Let \cref{ass:P} hold and consider a sequence $\{ (x^k,s^k) \}$ 
		of iterates generated by \cref{alg:ALM}.
		Then, each accumulation point $x^\ast$ of $\{ x^k \}$ 
		is feasible for \eqref{eq:P} if one of the following conditions holds:
		\begin{enumerate}
			\item\label{lem:ALM:bounded:mu}%
				$\{ \mu_k \}$ is bounded away from zero, or
			\item\label{lem:ALM:bounded:bound}%
				there exists some $B \in \R$ such that 
				$\LLslack_{\mu_k}(x^k, s^k,\hat{y}^k) \leq B$ for all $k\in\N$.
		\end{enumerate}
		In both situations, $(x^\ast,c(x^\ast))$ is an accumulation point
		of $\{ (x^k,s^k) \}$ which is feasible to \eqref{eq:ALMz:P}.
	\end{prop}
\end{mybox}
\begin{proof}
	Let $x^\ast \in \XX$ be an arbitrary accumulation point of $\{ x^k \}$ and 
	$\{x^k\}_K$ a subsequence such that $x^k \to_K x^\ast$.
	We need to show $c(x^\ast) \in D$
	under two circumstances.
	\begin{enumerate}
		\item If $\{ \mu_k \}$ is bounded away from zero, the conditions at
		\cref{state:ALM:if,state:ALM:iffailed} of \cref{alg:ALM} imply 
		that $\|c(x^k) - s^k\| \to 0$ for $k \to \infty$.
		By the upper bound $\|c(x^k) - s^k\| \geq \dist_D( c(x^k) )$ for all $k\in\N$, 
		due to $s^k \in D$,
		taking the limit $k \to_K \infty$ and continuity yield $\dist_D ( c(x^\ast) ) = 0$,
		hence $c(x^\ast) \in D$, i.e., $x^\ast$ is feasible to \eqref{eq:P}.
		Further, $s^k\to_K c(x^*)$ holds.
		\item 
		In case where $\{\mu_k\}$ is bounded away from zero, we can rely on the already proven
		first statement. 
		Thus, let us assume that $\mu_k\to 0$.
		By assumption, we have
		\begin{equation}\label{eq:intermediate_estimate_boundedness_of_AL_function}
			B
			\geq 
			\LLslack_{\mu_k}(x^k, s^k,\hat{y}^k) 
			= 
			q(x^k)
			+ 
			\frac{1}{2\mu_k} 
				\|c(x^k)+\mu_k\hat{y}^k-s^k\|^2
			- 
			\frac{\mu_k}{2}\|\hat{y}^k\|^2
		\end{equation}
		and $s^k\in D$ for all $k\in\N$.
		Rearranging terms yields the inequality
		\begin{equation*}
			q(x^k)
			+ 
			\frac{1}{2\mu_k} 
				\|c(x^k)+\mu_k\hat{y}^k-s^k\|^2
			\leq 
			B
			+ 
			\frac{\mu_k}{2}\|\hat{y}^k\|^2
		\end{equation*}
		for all $k\in\N$.
		Taking the lower limit $k \to_K \infty$ while respecting that $q$ is lsc and 
		$\{\hat y^k\}$ is bounded gives $x^\ast\in\dom q$.
		Particularly, $\{q(x^k)\}_{K}$ is bounded from below.
		Rearranging \eqref{eq:intermediate_estimate_boundedness_of_AL_function} yields
		\[
			\|c(x^k)+\mu_k\hat{y}^k-s^k\|^2
			\leq 
			2\mu_k(B-q(x^k))
			+ 
			\|\mu_k\hat{y}^k\|^2,
		\]
		and taking the upper limit $k\to_K\infty$ yields
		$\|c(x^k)-s^k\|\to_K 0$, again by boundedness of $\{\hat{y}^k\}$ and $\mu_k\to 0$.
		On the other hand, $c(x^k)\to_K c(x^\ast)$ follows by continuity, 
		and this gives $s^k\to_K c(x^\ast)$, since $D$ is closed and $s^k\in D$ for all $k\in\N$.
		Hence, $(x^\ast,c(x^\ast))$ is feasible to \eqref{eq:ALMz:P},
		i.e., $x^\ast$ is feasible to \eqref{eq:P}.
	\end{enumerate}
	The final statement of the proposition follows from the above arguments.
\end{proof}

The following convergence result provides fundamental theoretical support to \cref{alg:ALM}.
It shows that, under subsequential attentive convergence,
any feasible accumulation point is an AM-stationary point for \eqref{eq:P}.
\begin{mybox}
	\begin{thm}\label{thm:AMstationary}
		Let \cref{ass:P} hold and consider a sequence $\{ (x^k,s^k) \}$ of iterates 
		generated by \cref{alg:ALM} with $\varepsilon_k \to 0$.
		Let $(x^\ast,c(x^\ast))$ be an accumulation point  of $\{ (x^k,s^k) \}$ 
		feasible to \eqref{eq:ALMz:P} and
		$\{ (x^k,s^k)\}_K$ a subsequence such that $x^k \toattentive{q}_K x^\ast$ and $s^k \to_K c(x^\ast)$.
		Then, $x^\ast$ is an AM-stationary point for \eqref{eq:P}.
	\end{thm}
\end{mybox}
\begin{proof}
	Define $\zeta^k \coloneqq s^k - c(x^k)$ for all $k \in \N$.
	Then, from \cref{state:ALM:subproblem,state:ALM:y} of \cref{alg:ALM}, we have that
	\begin{align}
		- \nabla c(x^k)^\top y^k + \xi^k {}\in{}& \partial q(x^k) ,
		\label{thm:AMstationary:x} \\
		y^k + \nu^k {}\in{}& \limnormalcone_D(c(x^k) + \zeta^k)
		\label{thm:AMstationary:y}
	\end{align}
	for some $\xi^k \in \XX$, $\| \xi^k \| \leq \varepsilon_k$,
	and $\nu^k \in \YY$, $\| \nu^k \| \leq \varepsilon_k$; cf.\ \eqref{eq:epsstationary} and \eqref{eq:dLS}.
	Set $\lambda^k \coloneqq y^k+\nu^k$ and $\eta^k \coloneqq \nabla c(x^k)^\top\nu^k+\xi^k$ for all $k \in \N$.
	
	We claim that the four subsequences $\{x^k\}_K$, $\{\eta^k\}_K$, 
	$\{\lambda^k\}_K$
	and $\{\zeta^k\}_K$ satisfy the properties 
	in \cref{def:AMstationary} and therefore show that $x^\ast$ is an AM-stationary point 
	for \eqref{eq:P}.
	
	By construction, we have $x^k \toattentive{q}_K x^\ast$ 
	as well as $-\nabla c(x^k)^\top \lambda^k+\eta^k\in\partial q(x^k)$ and
	$\lambda^k\in\limnormalcone_D(c(x^k)+\zeta^k)$ for each $k \in \N$.
	Continuous differentiability of $c$ and $\|\xi^k\|,\|\nu^k\|\leq\varepsilon_k$ give $\|\eta^k\| \to_K 0$.
	Finally, $\zeta^k \to_K 0$ follows from $s^k \to_K c(x^\ast)$, $x^k \to_K x^\ast$ and continuity of $c$.
	
	Overall, this proves that $x^\ast$ is an AM-stationary point for \eqref{eq:P}.
\end{proof}

The additional assumption $x^k \toattentive{q}_K x^\ast$ in \cref{thm:AMstationary} is
trivially satisfied if $g$ is continuous on its domain since all iterates of
\cref{alg:ALM} belong to $\dom g$.
However, the following one-dimensional example illustrates how this additional requirement
appears to be indispensable in a discontinuous setting.
\begin{es}\label{ex:counterex}
	We consider $n\coloneqq m\coloneqq 1$ and set $D\coloneqq (-\infty,0]$,
	\[
		f(x)\coloneqq 0,
		\qquad
		g(x)\coloneqq 
		\begin{cases}
			x	&	\text{ if }x\leq 0,\\
			1-x	&	\text{ otherwise,}
		\end{cases}
		\qquad
		c(x)\coloneqq x.
	\]
	Note that $g$ is merely lsc at $x^\ast\coloneqq 0$,
	and that $\partial g(x^\ast)=[1,\infty)$; cf.\ \cref{fig:counterex_subdiff}.
	Although $x^\ast$ is the global maximizer of the associated	problem \eqref{eq:P},
	$x^\ast$ is not an M-stationary point.
	Since $\nabla f(x^\ast)=0$,
	$\nabla c(x^\ast)=1$
	and
	$\limnormalcone_D( c(x^\ast) )=\R_+$,
	there is no $y^\ast \in \limnormalcone_D( c(x^\ast) )$ such that $0 \in \nabla f(x^\ast) + \partial g(x^\ast) + \nabla c(x^\ast)^\top y^\ast$.
	Indeed, $x^\ast$ is not even AM-stationary.
	Possibly discarding early iterates,
	any sequence $\{x^k\}$ such that $x^k \toattentive{q} x^\ast$
	satisfies $x^k\leq 0$ for each $k\in\N$.
	Hence, we find $\partial q(x^k)\subset[1,\infty)$,
	$\nabla c(x^k)=1$ and 
	$\limnormalcone_D( c(x^k)+\zeta^k)\subset\R_+$ for
	each $\zeta^k\in\YY$ and $k\in\N$, showing that the distance between $0$
	and the set
	$\partial q(x^k)+\nabla c(x^k)^\top\limnormalcone_D( c(x^k)+\zeta^k)$
	is at least $1$.	
	
	We apply \cref{alg:ALM} with $\Ybounded\coloneqq \{0\}$, $\mu_0\coloneqq 1$,
	$\theta\coloneqq 1/4$ and $\kappa\coloneqq 1/2$.
	This may yield sequences $\{x^k\}$, $\{s^k\}$ and $\{\mu_k\}$ given by 
	$x^0 \coloneqq \mu_0$, $s^0:=0$, $x^k \coloneqq \mu_k \coloneqq 2^{1-k}$ and $s^k \coloneqq 0$
	for each $k\in\N$, $k\geq 1$; cf.\ \cref{fig:counterex_iterates}.
	Hence, we have $x^k \to x^\ast$ and,
	crucially, not $x^k \toattentive{q} x^\ast$.
\end{es}

\begin{figure}[tb!]
	\centering
	\begin{subfigure}[t]{0.48\linewidth}
		\centering
		\includetikz{counterex_subdiff}%
		\caption{Computation of $\partial g(0)$.}
		\label{fig:counterex_subdiff}
	\end{subfigure}
	\hfill
	\begin{subfigure}[t]{0.48\linewidth}
		\centering
		\includetikz{counterex_iterates}%
		\caption{Iterates $x^k$ for $k\in\{1,2,3\}$.}
		\label{fig:counterex_iterates}
	\end{subfigure}
	\caption{Visualizations for \cref{ex:counterex}.}
	\label{fig:counterex}
\end{figure}

The next result readily follows from \cref{cor:AMregular_points}
and \cref{thm:AMstationary}.
\begin{mybox}
	\begin{cor}
		Let \cref{ass:P} hold and consider a sequence $\{ (x^k,s^k) \}$ of iterates 
		generated by \cref{alg:ALM} with $\varepsilon_k \to 0$.
		Let $(x^\ast,c(x^\ast))$ be an accumulation point  of $\{ (x^k,s^k) \}$ 
		feasible to \eqref{eq:ALMz:P} and
		$\{ (x^k,s^k)\}_K$ a subsequence such that $x^k \toattentive{q}_K x^\ast$ and $s^k \to_K c(x^\ast)$.
		Furthermore, assume that $x^\ast$ is AM-regular for \eqref{eq:P}.
		Then, $x^\ast$ is an M-stationary point for \eqref{eq:P}.
	\end{cor}
\end{mybox}

We note that related results have been obtained in \cite[Thm~3.1]{chen2017augmented}
and \cite[Cor.~6.16]{KrugerMehlitz2021}.
In \cite{chen2017augmented}, however, the authors in most cases overlooked the
issue of attentive convergence in the definition of the limiting subdifferential for
discontinuous functions so that their findings are not reliable.

Constrained optimization algorithms aim at finding feasible points and minimizing 
the objective function subject to constraints.
Employing affordable local optimization techniques, 
one cannot expect to find global minimizers of any infeasibility measure.
Nevertheless, the next result proves that \cref{alg:ALM} with bounded 
$\{ \varepsilon_k \}$ finds stationary points of an infeasibility measure.
Notice that this property does not require $\varepsilon_k \to 0$, 
but only boundedness; cf. \cite[Thm~6.3]{birgin2014practical}.
\begin{mybox}
	\begin{prop} \label{lem:feasP}
		Let \cref{ass:P} hold and consider a sequence $\{ (x^k,s^k) \}$ of iterates 
		generated by \cref{alg:ALM} with $\{ \varepsilon_k \}$ bounded.
		Let $(x^\ast,s^\ast)$ be an accumulation point  of $\{ (x^k,s^k) \}$ 
		and $\{ (x^k,s^k)\}_K$ a subsequence such that $x^k \toattentive{q}_K x^\ast$ and $s^k \to_K s^\ast$.
		Then, $(x^\ast,q(x^\ast),s^*)$ is an M-stationary point of the feasibility problem
		\begin{equation}
			\minimize_{(x,\alpha,s) \in \epi q\times D}
			{}\quad{}
			\tfrac12\| c(x)-s \|^2.
			\label{eq:feasP}
		\end{equation}
		If $q$ is locally Lipschitz continuous
		at $x^\ast$, then $x^\ast$ is an M-stationary point of the constraint
		violation
		\begin{equation}
			\minimize_{(x,s)\in \XX\times D}
			{}\quad{}
			\tfrac12\| c(x)-s \|^2.
			\label{eq:constraint_violation}
		\end{equation}
	\end{prop}
\end{mybox}
\begin{proof}
	By \cref{lem:ALM:bounded:mu}, if $\{\mu_k\}$ is bounded away from zero, 
	$x^\ast$ is feasible for \eqref{eq:P} and $s^\ast=c(x^\ast)\in D$.
	Thus, $(x^\ast,q(x^\ast),c(x^\ast))$ is a global minimizer of \eqref{eq:feasP} 
	and $(x^\ast,c(x^\ast))$ is a global minimizer of \eqref{eq:constraint_violation}.
	By continuous differentiability of the objective function, M-stationarity with respect to
	both problems follows, see \cite[Prop.~5.1]{mordukhovich2006variational}.
	Hence, it remains to consider the case $\mu_k \to 0$.
	
	Owing to \cref{state:ALM:subproblem} of \cref{alg:ALM}, 
	for all $k \in \N$ it is 
	\begin{subequations}\label{eq:approx_stat_subproblem}
		\begin{align}
		\label{eq:approx_stat_subproblem_x}
		\xi^k 
		&\in 
		\partial q(x^k) 
		+ 
		\nabla c(x^k)^\top \left[ \hat{y}^k + (c(x^k) - s^k)/\mu_k \right],
		\\
		\label{eq:approx_stat_subproblem_s}
		\nu^k
		&\in 
		-\left[\hat{y}^k+(c(x^k) - s^k)/\mu_k\right]+\limnormalcone_D(s^k)
		\end{align}
	\end{subequations}
	for some $\xi^k \in \XX$, $\|\xi^k\| \leq \varepsilon_k$, and
	$\nu^k \in \YY$, $\|\nu^k\| \leq \varepsilon_k$; cf.\ \eqref{eq:dLS}.
	Particularly, \eqref{eq:approx_stat_subproblem_x} gives us
	\[
		(\xi^k-\nabla c(x^k)^\top[\hat y^k+(c(x^k)-s^k)/\mu_k],-1)
		\in 
		\limnormalcone_{\epi q}(x^k,q(x^k)).
	\]
	Multiplying by $\mu_k > 0$ and exploiting that $\limnormalcone_{\epi q}(x^k,q(x^k))$
	is a cone, we have
	\begin{equation}\label{eq:abstract_stationarity_feasibility_problem}
		(\mu_k\xi^k-\nabla c(x^k)^\top[c(x^k)+\mu_k\hat y^k-s^k],-\mu_k)
		\in 
		\limnormalcone_{\epi q}(x^k,q(x^k)).
	\end{equation}
	Furthermore, \eqref{eq:approx_stat_subproblem_s} yields
	\begin{equation}\label{eq:abstract_ststionarity_feasibility_problem_s}
		\mu_k(\nu^k+\hat{y}^k)
		+
		c(x^k)-s^k
		\in 
		\limnormalcone_D(s^k)
	\end{equation}
	since $\limnormalcone_D(s^k)$ is a cone.
	Taking the limit $k\to_K\infty$ in \eqref{eq:abstract_stationarity_feasibility_problem}
	and \eqref{eq:abstract_ststionarity_feasibility_problem_s},
	the robustness of the limiting normal cone, $x^k \toattentive{q}_K x^\ast$
	and boundedness of $\{\hat{y}^k\}$, $\{\xi^k\}$ and $\{\nu^k\}$ yield
	\begin{equation}\label{eq:stationarity_constraint_violation}
		\begin{aligned}
		(-\nabla c(x^\ast)^\top[c(x^\ast)-s^\ast],0)
		&\in 
		\limnormalcone_{\epi q}(x^\ast,q(x^\ast)),
		\\
		c(x^\ast)-s^\ast
		&\in 
		\limnormalcone_D(s^\ast).
		\end{aligned}
	\end{equation}
	Keeping the Cartesian product rule for the computation of limiting normals
	in mind, see \cite[Prop.~1.2]{mordukhovich2006variational},
	$(x^\ast,q(x^\ast),s^\ast)$ is an M-stationary point of
	\eqref{eq:feasP}.
	
	Finally, assume that $q$ is locally Lipschitz
	continuous at $x^\ast$. Then, due to
	\cite[Cor.~1.81]{mordukhovich2006variational}, we have
	\[
		(y^\ast,0)\in\limnormalcone_{\epi q}(x^\ast,q(x^\ast))
		\quad\Rightarrow\quad
		y^\ast=0,
	\]
	so that the above arguments already show M-stationarity of
	$(x^\ast,s^\ast)$ for \eqref{eq:constraint_violation}.
	\qedhere
\end{proof}

In case where $D$ is convex, the assertion of \cref{lem:feasP} can be slightly strengthened.
\begin{mybox}
	\begin{cor}
		Let $D$ be convex, let \cref{ass:P} hold and consider a sequence $\{ (x^k,s^k) \}$ of iterates 
		generated by \cref{alg:ALM} with $\{ \varepsilon_k \}$ bounded.
		Let $(x^\ast,s^\ast)$ be an accumulation point  of $\{ (x^k,s^k) \}$ 
		and $\{ (x^k,s^k)\}_K$ a subsequence such that $x^k \toattentive{q}_K x^\ast$ and $s^k \to_K s^\ast$.
		Then, $(x^\ast,q(x^\ast))$ is an M-stationary point of the feasibility problem
		\begin{equation*}
			\minimize_{(x,\alpha) \in \epi q}
			{}\quad{}
			\tfrac12 \dist_D^2(c(x)).
		\end{equation*}
		If $q$ is locally Lipschitz continuous
		at $x^\ast$, then $x^\ast$ is an M-stationary point of the constraint
		violation
		\begin{equation*}
			\minimize_{x\in \XX}
			{}\quad{}
			\tfrac12 \dist_D^2(c(x)).
		\end{equation*}
	\end{cor}
\end{mybox}
\begin{proof}
	We proceed as in the proof of \cref{lem:feasP} in order to come up with
	\eqref{eq:stationarity_constraint_violation}.
	By convexity of $D$, $c(x^\ast)-s^\ast\in\limnormalcone_D(s^\ast)$ is equivalent
	to $s^*\in\proj_D(c(x^\ast))$.
	Thus, the assertion follows from \cref{lem:squared_distance}.
\end{proof}

	\subsection{Termination Criteria}\label{sec:ALM:termination}
	\Cref{state:ALM:subproblem} involves the minimization of the augmented Lagrangian function defined in \eqref{eq:L}.
Then, the dual update at \cref{state:ALM:y} allows to draw conclusions with respect to the original problem \eqref{eq:P}, as 
\cref{thm:AMstationary} shows that accumulation points of sequences generated 
by \cref{alg:ALM} are AM-stationary under mild assumptions.

Owing to \eqref{thm:AMstationary:x}--\eqref{thm:AMstationary:y} and recalling the AM-stationarity conditions \eqref{eq:AMstationary}, one may select a null sequence $\{ \varepsilon^k \} \subset \R_{++}$ at \cref{state:ALM:ysafe}.
Then, given some user-defined tolerances $\varepsilon^{\text{dual}}, \varepsilon^{\text{prim}} > 0$, it is reasonable to declare successful convergence when the conditions
\begin{equation*}
	\varepsilon^k \leq \varepsilon^{\text{dual}}
	\quad\text{and}\quad
	\| c(x^k) - s^k \| \leq \varepsilon^{\text{prim}}
\end{equation*}
are satisfied.
\Cref{thm:AMstationary} demonstrates that these termination criteria (the latter, in particular) are satisfied in finitely many iterations 
if any subsequence of $\{(x^k,s^k)\}$ accumulates at a feasible point $(x^\ast,c(x^\ast))$ 
of \eqref{eq:ALMz:P}.
As this might not be the case, a mechanism for (local) infeasibility detection is needed, and usually included in practical implementations; see \cite{armand2019rapid,burke2014sequential}.

Given some tolerances, \cref{alg:ALM} can be equipped with relaxed conditions on decrease requirements at \cref{state:ALM:if} and optimality at \cref{state:ALM:subproblem}.
At \cref{state:ALM:ysafe} the inner tolerance $\varepsilon^k$ can stay bounded away from zero, as long as $\varepsilon^k \leq \varepsilon^{\text{dual}}$ for large $k \in \N$.
Similarly, the condition at \cref{state:ALM:if} can be relaxed by adding the (inclusive) possibility that $\| c(x^k) - s^k \| \leq \varepsilon^{\text{prim}}$.
Finally, at \cref{state:ALM:ifpassed} a nonmonotone update is allowed, namely the penalty parameter can be increased, as long as some watchdog procedures are in place to avoid cycling \cite{birgin2012augmented}.

	\subsection{Inner Problem and Solver}\label{sec:inner}
	In this section we elaborate upon \cref{state:ALM:subproblem} of \cref{alg:ALM} that aims at minimizing the augmented Lagrangian function $\LLslack_\mu(\cdot,\cdot,y)$ defined in \eqref{eq:Lz}.
To this end, let us take a closer look at the structure of this 
subproblem.

Using the decomposition $\LLslack_\mu(\cdot,\cdot,y) = f^{\textup{S}}(\cdot,\cdot) + g^{\textup{S}}(\cdot,\cdot)$ with component functions $f^{\textup{S}}\colon\XX\times\YY\to\R$
and $g^{\textup{S}}\colon\XX\times\YY\to\Rinf$ given by
\begin{align}
	f^{\textup{S}}(x,s)
	{}\coloneqq{}&
	f(x) + \frac{1}{2\mu} \| c(x) + \mu y - s \|^2 - \frac{\mu}{2}\|y\|^2 ,
	\label{eq:fz} \\
	g^{\textup{S}}(x,s)
	{}\coloneqq{}&
	g(x) + \indicator_D(s) ,
	\label{eq:gz}
\end{align}
one immediately sees that this split recovers the classical setting of an 
\emph{unconstrained} composite optimization problem with $f^{\textup{S}}$ being
continuously differentiable, while $g^{\textup{S}}$ is merely lsc, but of
a particular structure. In principle, proximal gradient-type methods can therefore
be applied as approximate solvers for our subproblems, see \cite{beck2017first} for an
introduction of this class of methods. A standing assumption of the corresponding convergence theory in \cite{beck2017first} and all previous works on 
proximal gradient-type methods,
however, is a global Lipschitz condition regarding the gradient of the
smooth part $f^{\textup{S}}$. Note that this gradient is given by
\begin{equation*}
	\nabla f^{\textup{S}}(x,s) = \begin{bmatrix}
		\nabla f(x) + \frac{1}{\mu} \nabla c(x)^\top \left[ c(x) + \mu y - s \right] \\
		- \frac{1}{\mu} \left[ c(x) + \mu y - s \right] 
	\end{bmatrix} .
\end{equation*}
Observe that our standing assumptions from \cref{ass:P} imply that this gradient
is locally Lipschitz continuous, but they do not guarantee global Lipschitzness 
in general. Fortunately, some recent contributions on
proximal gradient-type methods show that these methods also work under
suitable assumptions if the smooth term has a locally Lipschitz gradient
only; cf.\ \cite{bauschke2017descent,demarchi2022proximal,KanzowMehlitz2022} for more details. Consequently, these proximal gradient-type methods offer a viable way to solve the augmented Lagrangian subproblems, even for fully nonconvex problems.
Let us also mention that, at least in \cite{demarchi2022proximal,KanzowMehlitz2022}, it has
been verified that accumulation points of sequences generated by
proximal gradient-type methods are stationary while along the associated subsequence,
the iterates are $\varepsilon_k$-stationary for a null sequence $\{\varepsilon_k\}$.
This requirement is essential in \cref{alg:ALM}.

For a practical implementation of these proximal methods, it is advantageous to exploit 
the particular structure of the nonsmooth term $g^{\textup{S}}$. In fact, due to the
separability of $g^{\textup{S}}$ with respect to $x$ and $s$, it follows that
the corresponding proximal mapping is easily computable. More precisely, one obtains
\begin{equation*}
	\prox_{\gamma g^{\textup{S}}}(x,s) 
	= 
	\begin{bmatrix}
		\prox_{\gamma g}(x) \\
		\proj_D(s)
	\end{bmatrix}
\end{equation*}
for any $\gamma \in (0,\gamma_g)$.

Though the proximal-type approach is used in our numerical setting (see the
next section for some more details), we stress that there exist other 
candidates for the numerical solution of the resulting augmented Lagrangian
subproblems. To this end, recall that the previous discussion looked at
these subproblems as an unconstrained composite optimization problem.
Alternatively, we may view these subproblems from the point of view of
machine learning, where (essentially) the same class of optimization problems is
solved by (possibly) different techniques. We refer the interested
reader to \cite{sra2011optimization,wright2022optimization}
for a survey of optimization methods
for machine learning and data analysis problems.
These techniques might be applicable
very successfully at least in certain situations. For example, if the
smooth term $f^{\textup{S}}$ is convex (the gradient does not have to
be globally Lipschitz), whereas the nonsmooth term $g^{\textup{S}}$ is still
just assumed to be lsc (and not necessarily convex), it 
is possible to adapt the idea of cutting plane methods to this setting by
applying the cutting plane technique to $f^{\textup{S}}$ only, whereas one
does not change the nonsmooth term. The resulting subproblems then use
a piecewise affine lower bound for the function $f^{\textup{S}}$ and add the
(possibly complicated) function $g^{\textup{S}}$. Of course, and similar
to the proximal gradient-type approaches, these subproblems need to be
easily solvable for the overall augmented Lagrangian method to be
efficient, and this, in general, is true only for particular classes of
problems; cf.\ \cref{sec:numresults}.

\section{Numerical Examples}\label{sec:numresults}
	This section presents a numerical implementation of \cref{alg:ALM} and discusses its behavior on some illustrative examples, showcasing the flexibility offered by the constrained composite programming framework.
In particular, we consider challenging problems where the cost function is nonsmooth and
nonconvex or where the constraints are inherently nonconvex by a disjunctive structure
of the respective set $D$.
In \cref{sec:numresults:rosenbrock} we demonstrate the benefit of accelerated proximal-gradient methods for solving the subproblems by means of a simple two-dimensional problem where
a nonsmooth variant of the Rosenbrock function is minimized over a set of combinatorial
structure.
Next, \cref{sec:numresults:scsto} is dedicated to a binary optimal control problem with nonlinear dynamics, free final time and switching costs, where we display and discuss weaknesses of our approach.
\Cref{sec:numresults:sparseportfolio} deals with a test collection of portfolio optimization
problems from \cite{frangioni2007sdp} which are equipped with a nonconvex 
sparsity-promoting term in the objective function.
Finally, in \cref{sec:numresults:matrixcompletionminrank} 
we address a class of matrix recovery problems discussed e.g.\ in
\cite{shen2018penalty} where the rank of the unknown matrix has to be minimized.

	\subsection{Implementation}\label{sec:numresults:implementation}
	We have implemented the proposed Augmented Lagrangian Solver (\alps{}) as part of an open-source software package in the Julia language \cite{bezanson2017julia}.
\alps{} can solve constrained composite problems of the form \eqref{eq:P} and is available online at
\begin{center}
	\BazingaFullLink{},
\end{center}
together with the examples presented in the following sections.
\alps{} can be used to solve, in the sense of \cref{sec:ALM:termination}, a wide spectrum of optimization problems, requiring only first-order primitives, \ie, gradient, proximal mapping and projections.
By default, \alps{} invokes \panoc{+} \cite{demarchi2022proximal} for solving the augmented Lagrangian subproblems at \cref{state:ALM:subproblem} of \cref{alg:ALM}, possibly inexactly and up to stationarity, using the implementation offered by \ProximalAlgorithmsLink{} \cite{stella2017proximalalgorithms};
see \cref{sec:appendix} for more details.
The method is implemented matrix-free, that is, the constraint Jacobian $\nabla c$ does not need to be explicitly formed as only Jacobian-vector products $\nabla c(x)^\top v$ are required.

The solver requires the data functions $f$, $g$, $c$ and constraint set $D$ specified as objects returning the oracles discussed at the end of \cref{sec:introduction}.
Further, the initialization requires a primal-dual starting point $(x^{\text{init}},y^{\text{init}}) \in \XX \times \YY$.
The default safeguarding set $\Ybounded$ in $\YY$ is $\Ybounded = [-y_{\max}, y_{\max}]^m$, with $y_{\max} = 10^{20}$, and the safeguarded dual estimate $\hat{y}^k$ at \cref{state:ALM:ysafe} is chosen as the projection of $y^{k-1}$ onto $\Ybounded$; of $y^{\text{init}}$ for $k=0$.
User override of this oracle allows for tailored choices of $\Ybounded$, possibly exploiting the structure of $D$ \cite{sopasakis2020open}.

\alps{} initializes \cref{alg:ALM} by overwriting $x^{\text{init}}$ with an arbitrary element of $\prox_{\gamma g}(x^{\text{init}}) \subset \dom q$, where $\gamma = \epsilon_M$ and $\epsilon_M$ denotes the machine epsilon of a given floating-point system.
The examples presented in the following are in double precision (Float64), so $\epsilon_M \approx 2.22 \cdot 10^{-16}$.
The inner tolerances $\varepsilon_k$ at \cref{state:ALM:ysafe} are constructed as a sequence of decreasing values, defined by the recurrence
\begin{equation*}
	\varepsilon_{k+1} = \max \{ \kappa_\varepsilon \varepsilon_k, \varepsilon^{\text{dual}} \} ,
\end{equation*}
starting from $\varepsilon_0 \coloneqq (\varepsilon^{\text{dual}})^{\frac{1}{3}}$ and given some $\varepsilon^{\text{dual}}, \kappa_\varepsilon \in (0,1)$ \cite{birgin2012augmented}.
The initial penalty parameter $\mu_0$ is automatically chosen by default, similarly to \cite[Eq. 12.1]{birgin2014practical}.
Given $x^{\text{init}} \in \dom q$, we evaluate the constraints $c^{\text{init}} \coloneqq c(x^{\text{init}})$, select an arbitrary element $s^{\text{init}} \in \proj_D( c^{\text{init}} )$ and compute the vector $\Delta^{\text{init}} \coloneqq c^{\text{init}} - s^{\text{init}}$.
Then, the vector $\mu_0 \in \YY$ of penalty parameters is selected componentwise as follows:
\begin{equation*}
	(\mu_0)_i \coloneqq \max\left\{ 10^{-8}, \min\left\{ \frac{1}{10} \frac{\max\{ 1, (\Delta_i^{\text{init}})^2/2 \}}{\max\{ 1, q(x^{\text{init}}) \}}, 10^8 \right\} \right\} ,
\end{equation*}
effectively scaling the contribution of each constraint \cite{birgin2014practical,conn1991globally}.
Then, according to the overall feasibility-complementarity of the iterate, the penalty parameters are updated in unison at \cref{state:ALM:iffailed}, since using a different penalty parameter for each constraint is theoretically worse than using a common parameter \cite[\S 3.4]{andreani2021best}; we set $\mu_{k+1} \coloneqq \kappa_\mu \mu_k$, for some fixed $\kappa_\mu \in (0,1)$.
At the $k$th iteration, the subsolver at \cref{state:ALM:subproblem} is warm-started
from the previous estimate $(x^{k-1},s^{k-1}) \in \dom q \times D$;
from $(x^{\text{init}}, s^{\text{init}})$ for $k=0$.

The default parameters in \alps{} are $\theta = 0.8$, $\kappa_\mu = 0.5$ and $\kappa_\varepsilon = 0.1$, termination tolerances $\varepsilon^{\text{prim}} = \varepsilon^{\text{dual}} = 10^{-6}$ and a maximum number of (outer) iterations, whose default value is $100$.

	\subsection{Nonsmooth Rosenbrock and Either-Or Constraints}\label{sec:numresults:rosenbrock}
\begin{figure}[tb!]
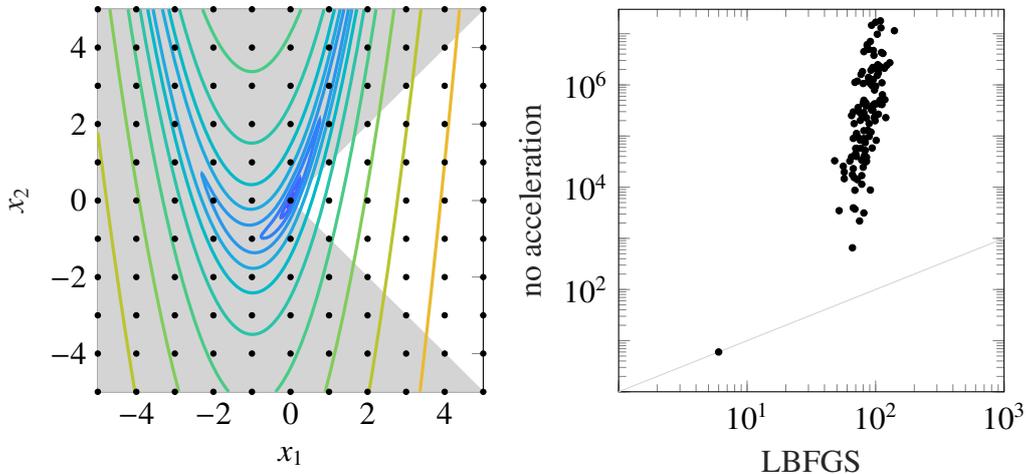

	\centering
		\includetikz{rosenbrock}%
	\caption{%
		Setup and results for the illustrative problem \eqref{eq:rosenbrock}.
		Left: Feasible region (gray background), objective contour lines, global minimizer $x^\ast = (0,0)$ and grid of starting points.
		Right: Comparison of inner iterations needed without acceleration against LBFGS acceleration; each mark corresponds to a starting point and the gray line has unitary slope.%
	}%
	\label{fig:rosenbrock}
\end{figure}

	Let us consider a two-dimensional optimization problem involving a non\-smooth Rosenbrock-like objective function and either-or constraints, namely set-membership constraints entailing an inclusive disjunction.
It reads
\begin{equation}\label{eq:rosenbrock}
	\minimize_{x}
	~{}
	10 (x_2 + 1 - (x_1 + 1)^2)^2 + |x_1|
	\quad{}
	\stt
	~{}
	x_2 \leq - x_1
	{}~\vee~{}
	x_2 \geq x_1
\end{equation}
and admits a unique (global) minimizer $x^\ast = (0, 0)$.
The feasible set is nonconvex and connected; see \cref{fig:rosenbrock}.
We cast \eqref{eq:rosenbrock} into the form of \eqref{eq:P} by defining the data functions as
\begin{equation*}
	f(x) \coloneqq 10 (x_2 + 1 - (x_1 + 1)^2)^2 ,
	{}\qquad{}
	g(x) \coloneqq |x_1| ,
	{}\qquad{}
	c(x) \coloneqq
	\begin{pmatrix}
		-x_1-x_2 \\
		-x_1+x_2
	\end{pmatrix} ,
\end{equation*}
and let the constraint set be $D \coloneqq \Deo$, where the (nonconvex) set
\begin{equation*}
	\Deo
	{}\coloneqq{}
	\{ (a,b) \suchthat a \geq 0 \vee b \geq 0 \}
	{}={}
	\{(a,b) \suchthat a \geq 0\} \cup \{(a,b) \suchthat b \geq 0\}
\end{equation*}
describes the either-or constraint.

We consider a uniform grid of $11^2=121$ starting points $x^0$ in $[-5,5]^2$ 
and let the initial dual estimate be $y^0 = 0$.
Also, we compare the performance of \alps{} by solving the subproblems using \panoc{+} without or with (LBFGS) acceleration; see the last paragraph of \cref{sec:appendix} for more details.

\alps{} solves all the problem instances, approximately (tolerance $10^{-3}$ in Euclidean distance) reaching $x^\ast = (0,0)$ in all cases.
\cref{fig:rosenbrock} depicts the feasible region of \eqref{eq:rosenbrock}, some contour lines of its objective function and the grid of starting points $x^0$.
Over all problems, \alps{} with no acceleration takes at most $17 \, 870 \, 346$ (cumulative) inner iterations to find a solution (median $291 \, 756$), whereas with LBFGS directions only $140$ inner iterations are needed at most (median $86$).
A closer look at \cref{fig:rosenbrock} indicates that not only the accelerated method usually requires far less iterations, but also that its behavior is more consistent, as the majority of cases spread over a narrow interval.
These results support the claim that (quasi-Newton) acceleration techniques can give a mean to cope with bad scaling and ill-conditioning \cite{stella2017simple,themelis2018proximal}.

	\subsection{Sparse Switching Time Optimization}\label{sec:numresults:scsto}
	Constrained composite programming offers a flexible language for modeling a variety of problems.
In this section we consider the sparse binary optimal control of Lotka-Volterra dynamics.
Known as the fishing problem \cite[\S 6.4]{sager2005numerical}, it is typically stated as
\begin{align}
	\minimize_{x,u}
	\quad{}&
	\int_0^T \| x(t) - 1 \|^2 \mathrm{d}t \label{eq:scsto:fishing}\\
	\stt
	\quad{}&
	\dot{x}_1(t) {}={} x_1(t) [ - c_1 u(t) - x_2(t) + 1 ] &&\text{for a.e.}\; t\in[0,T], \nonumber\\
	&
	\dot{x}_2(t) {}={} x_2(t) [ - c_2 u(t) + x_1(t) - 1 ] &&\text{for a.e.}\; t\in[0,T], \nonumber\\
	&
	x(0) = x_0, \nonumber\\
	&
	u(t) \in \{0,1\} &&\text{for}\; t\in[0,T], \nonumber
\end{align}
where final time $T = 12$, initial state $x_0 \coloneqq (0.5, 0.7)$ and parameters $c_1 = 0.4$, $c_2 = 0.2$ are given and fixed.
In order to showcase the peculiar features of \eqref{eq:P}, we focus on a variant of the fishing problem with switch costs and free, although constrained, final time.
First, the problem is reformulated as a finite-dimensional one by adopting the switching time optimization approach, that consists in optimizing the times at which the control input changes, given a fixed sequence of $N$ admissible controls \cite[\S 5.2]{sager2005numerical}.
We call switching intervals the time between these switching times and collect them in a vector $\tau \in \R^N$.
Clearly, they must take nonnegative values and sum up to the final time $T$.
Furthermore, considering the chattering solution exhibited by the fishing problem \cite[\S 6.5]{sager2005numerical}, we introduce switch costs to penalize solutions that show frequent switching of the binary control trajectory, yielding more practical results.
Following \cite{demarchi2020constrained}, \cite[Ch.~2]{demarchi2021dissertation}, switch costs can be interpreted as a regularization term and modeled using the $\ell_0$ quasi-norm of the switching intervals, effectively counting how many control inputs in the given control sequence are active.
The resulting problem formulation reads
\begin{equation}\label{eq:scsto}
	\minimize_{\tau}
	\quad{}
	f(\tau) + \indicator_{\R_+^N}(\tau) + \sigma \| \tau \|_0
	\qquad{}
	\stt
	\quad{}
	1_N^\top \tau \in D . 
\end{equation}
Here, the smooth cost function $f$ returns the tracking cost, by integrating the dynamics, starting from the initial state, for the given sequence of control inputs and switching intervals.
The nonnegativity constraint $\indicator_{\R_+^N}$ and sparsity-promoting cost $\sigma \|\cdot\|_0$ form the nonsmooth cost function $g$ in \eqref{eq:P}; despite $g$ being nonconvex and discontinuous, its proximal mapping can be easily evaluated \cite[\S 3.2]{demarchi2020constrained}.
The nonnegative parameter $\sigma$ controls the impact of the $\ell_0$ regularization and can be interpreted as the switching cost.
The only constraint remained explicit is the one on the final time $T \coloneqq 1_N^\top \tau$.
Hence, the constraint set $D \subset \R_+$ is constituted by the admissible values for $T$.

We consider the binary control sequence $\{ 0,1,0,\dots,1 \}$ with $N \coloneqq 24$ intervals.
A background time grid with $n = 200$ points is adopted to integrate dynamics and evaluate sensitivities, following the linearization approach of \cite{stellato2017second}.
We solve \eqref{eq:scsto} for increasing values of the switching cost parameter $\sigma \in \{ 10^{-6}, 10^{-5}, \dots, 10 \}$.
For the first problem, the initial guess $\tau^0$ corresponds to uniform switching intervals with the final time $T = 12$ usually fixed in \eqref{eq:scsto:fishing}.
Then, following a continuation approach, a solution is adopted as initial guess for the subsequent problem, but always with dual estimate $y^0 = 0$.
Moreover, we consider two cases for the constraint set $D$.
First, we let $D \coloneqq [0,15]$ and \alps{} returns solutions whose final time reaches values around $T \approx 12$.
Then, we consider a second case with the disconnected constraint set $D \coloneqq [5,10] \cup [13,15]$, so to impact on the solution; in this case the returned final times are $T \approx 13$.

\begin{figure}[tb]
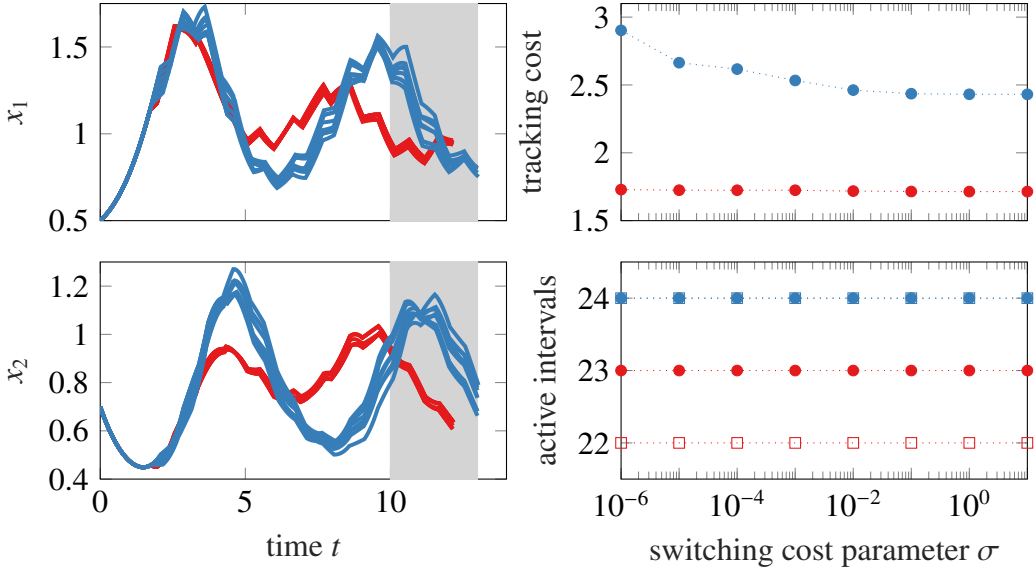

	\centering
	\includetikz{scsto}%
	\caption{%
		Results for the illustrative problem \eqref{eq:scsto} using switching time optimization with a sequence of $24$ binary controls and several values for the switching cost parameter $\sigma$.
		Left: Prohibited region for the final time (gray background) and state trajectories with (blue) or without (red) constraint.
		Right: Comparison of the resulting tracking cost and number of nonzero variables, corresponding to active intervals (circle).
		Identical control trajectories can be obtained with fewer active intervals (square), yielding lower switching cost.%
	}%
	\label{fig:scsto}
\end{figure}

\alps{} is able to find reasonable solutions that satisfy the constraints, despite the nonconvexity of the switching time approach \cite[Apx B.4]{sager2005numerical}, the discrete nature of the sparse regularizer and the constraint set $D$ being disconnected.
It should be stressed, however, that there are no guarantees on the quality of these solutions and, in fact, the solutions found by \alps{} are poor in terms of objective value, as we are about to show.

The state trajectories are depicted in \cref{fig:scsto}, for both cases, along with a comparison of the tracking cost and number of active intervals against the switching cost parameter $\sigma$.
First, we observe that the trajectories are not strongly affected, despite the dramatic increase of $\sigma$ (relative to the tracking cost).
Moreover, the solver performs only few iterations, needed to adjust the dual estimate and verify the termination criteria.
In practice, the iterates remain trapped around a minimizer with high objective value, and a huge value of $\sigma$ is required for jumping to a lower objective value.
This becomes apparent looking at $\| \tau \|_0$, namely the number of active intervals.
Given a sequence of control inputs, several choices of switching intervals can give the same state trajectory, hence the same tracking cost.
Among these, we would expect the solver to return one with minimum number of nonzeros.
For instance, vectors of switching intervals in the form $(\alpha+\beta,0,0,\dots)$ and $(0,0,\alpha+\beta,\dots)$ should be preferred over $(\alpha,0,\beta,\dots)$, for they yield the same control trajectory whilst having fewer nonzero elements.
The solutions returned by \alps{} are compared against equivalent although sparser ones in \cref{fig:scsto}.
Clearly, and not surprisingly, the solutions obtained are far from being globally optimal.

	\subsection{Sparse Portfolio Optimization}\label{sec:numresults:sparseportfolio}
	Let us consider portfolio optimization problems in the form
\begin{equation}\label{eq:sparseportfolio}
	\begin{aligned}
		\minimize_x
		\quad{}{}&
		\frac{1}{2} x^\top Q x + \alpha \|x\|_0 \\
		\stt
		\quad{}{}&
		\mu^\top x \geq \varrho
		,\quad
		1_n^\top x = 1
		,\quad
		0 \leq x \leq u .
	\end{aligned}
\end{equation}
The problem data $Q \in \R^{n \times n}$ and $\mu \in \R^n$ 
denote the covariance matrix and the mean of $n\in\N$ possible assets, respectively, 
while $\varrho \in \R$ is a lower bound for the expected return.
Furthermore, $u \in \R^n$ provides an upper bound for the individual assets within the portfolio.
Aiming at a sparse portfolio, 
and in contrast with cardinality-constrained formulations, see e.g.\ \cite{jia2021augmented}, 
we use the $\ell_0$ quasi-norm as a regularization term that penalizes the 
number of chosen assets within the portfolio.

We reformulate the model in the form of \eqref{eq:P} by letting 
$f$ be the quadratic cost, $g$ the nonsmooth cost and indicator of the bounds, 
$\func{c}{\R^n}{\R^m}$, $m \coloneqq 2$, defined by $c(x) \coloneqq [\mu, 1_n]^\top x$ 
and $D \coloneqq [\varrho, \infty) \times \{1\}$.

Through a mixed-integer quadratic program formulation of \eqref{eq:sparseportfolio},
which can be obtained via the theory provided in \cite{FengMitchellPangShenWaechter2018},
we compute a solution using CPLEX \cite{cplex2009v12}, for comparison.
Based on our experiences from \cref{sec:numresults:scsto},
we also solve \eqref{eq:sparseportfolio} using a continuation procedure: 
the $\ell_0$ minimization is warm-started at a primal-dual point found replacing 
the discontinuous $\ell_0$ function with either the norm $\ell_1\coloneqq\|\cdot\|_1$ or 
the $p$-th power of the $\ell_p$ quasi-norm, \ie,
$\ell_p^p\coloneqq\|\cdot\|_p^p$ ($p = 0.5$) and solving the corresponding problem.
Notice that \eqref{eq:sparseportfolio} with the $\ell_0$- replaced by
the $\ell_1$-term boils down to a convex quadratic program; in fact, it is $\|x\|_1 = 1$ 
for each feasible point of \eqref{eq:sparseportfolio} 
by the nonnegativity and equality constraints.

The data $Q$, $\mu$, $\varrho$ and $u$ is taken from the test problem collection
\cite{frangioni2007sdp}, which has been created randomly and is available online \cite{frangioni2021mean}.
Here, we used all 30 test instances of dimension $n \coloneqq 200$ 
and the two different values $\alpha \in \{10, 100\}$ for each problem.

\begin{figure}[tb]
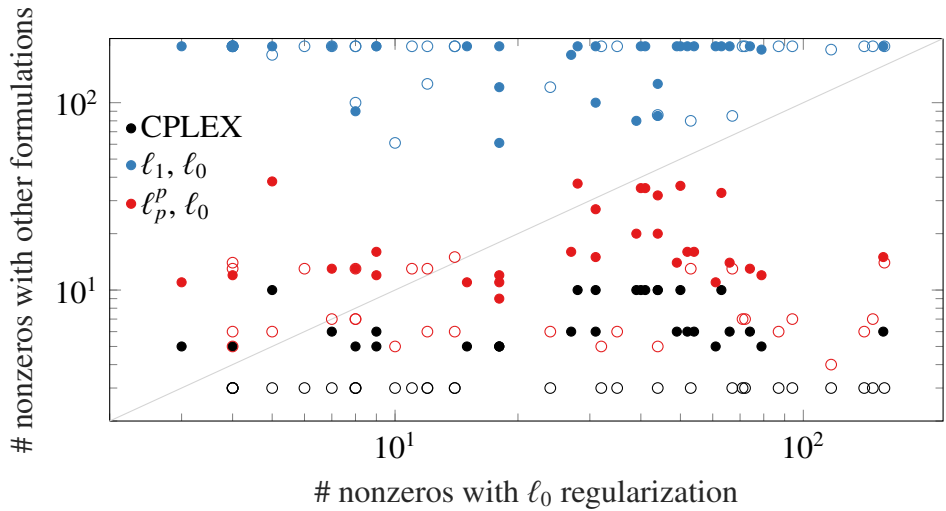

	\centering
	\includetikz{portfolio}%
	\caption{%
		Results for the portfolio problem \eqref{eq:sparseportfolio}:
		Comparison of the solutions found with $\ell_0$ regularization against those obtained with CPLEX
		and $\ell_0$ warm-started with $\ell_1$ or $\ell_p^p$, with $p=0.5$.
		We depict the number of nonzero entries of the solutions 
		returned for $\alpha=10$ (dot) and $\alpha=100$ (circle).
		The gray line has unitary slope.%
	}%
	\label{fig:portfolio}%
\end{figure}

The results of our experiments are depicted in \cref{fig:portfolio}.
Let us mention that \alps{} solved all problem instances, in the sense that it returned primal-dual pairs satisfying the termination criteria of \cref{sec:ALM:termination}.
Below, we comment on some median values for our experiments with
parameters $\alpha=10/100$:
a direct use of $\ell_0$ minimization resulted in
$10/13$ outer and $908/1633$ inner iterations,
while warm-starting with the continuous $\ell_p^p$ function
required
$13/9$ outer and $686/1830$ inner iterations.
Let us point the reader's attention to the fact that
the $\ell_p^p$-warm-started $\ell_0$ minimization 
did not affect the solution sparsity, \ie, the numbers of nonzero
components of the obtained solutions were the same
with and without an additional round of $\ell_0$ minimization after the
$\ell_p^p$ warm-start.
Although one cannot expect to find a global minimum in general, 
we recall that the standard $\ell_1$ regularization does not work in this example,
as confirmed by the poor performance depicted in \cref{fig:portfolio},
whereas the nonconvex $\ell_p^p$ penalty already leads to very sparse solutions.

	\subsection{Matrix Completion with Minimum Rank}\label{sec:numresults:matrixcompletionminrank}
	For some $\ell\in\N$, $\ell\geq 2$, 
let us consider $N \in \N$ points $x_1,\ldots,x_N \in \R^\ell$
and define a block matrix $X \in \R^{N \times \ell}$ by means of
$X \coloneqq [x_1, x_2, \dots, x_N]^\top$.
Let $\Delta \in \R^{N \times N}$ denote the Euclidean distance matrix associated with these points,
given by $\Delta_{ij} \coloneqq \|x_i - x_j\|^2 = (x_i-x_j)^\top(x_i-x_j)$
for all $i, j \in \mathcal{I} \coloneqq \{1,\dots,N\}$.
We aim at recovering $X$ based on a partial knowledge of $\Delta$.
In particular, we assume that $\Omega \subset \mathcal{I}^2$ is a set
of pairs such that only the entries $\Delta_{ij}$, $(i,j)\in\Omega$, 
of $\Delta$ are known.

Following \cite{shen2018penalty}, we lift the problem by introducing 
a symmetric matrix $B \coloneqq X X^\top$ 
whose rank is, by construction, smaller than or equal to $\ell$.
Hence, we seek a matrix $B \in \R^{N \times N}$ that satisfies the symmetry constraint 
$B = B^\top$ and the distance constraints associated with the observations,
i.e., $B_{ii}+B_{jj}-B_{ij}-B_{ji}=\Delta_{ij}$ has to hold for all $(i,j)\in\Omega$.
Among these admissible matrices, those with minimum rank are preferred.

Let us consider problems of type
\begin{equation}\label{eq:matrixcompletionminrank}
	\begin{aligned}
	\minimize_{B}
	\quad{}{}&
	g(B) \\
	\stt
	\quad{}{}&
	B_{ii} + B_{jj} - B_{ij} - B_{ji} = \Delta_{ij} & & \forall (i,j) \in \Omega,\\
	&
	B_{ij} = B_{ji} & & \forall i, j \in \mathcal{I}, j < i 
	\end{aligned}
\end{equation}
where the function $\func{g}{\R^{N \times N}}{\R}$ encodes a matrix regularization term.
In the following, we consider $g \coloneqq \rank \coloneqq \| \sigma(\cdot) \|_0$, 
the nuclear norm $g \coloneqq \|\cdot\|_\ast \coloneqq \sum_i \sigma_i(\cdot)$ 
or the $p$-powered Schatten $p$-quasi-norm 
$g \coloneqq \|\cdot\|_p^p \coloneqq \sum_i \sigma_i(\cdot)^p$, $p\in(0,1)$,
where $\sigma(A)$ denotes the vector of singular values of a matrix $A$.
In our experiments rank and singular values are numerically evaluated using Julia's LinearAlgebra functions \texttt{rank} and \texttt{svd}, respectively.
Notice in particular that the rank of a matrix $A$ is computed by counting how many singular values of $A$ have magnitude greater than a numerical tolerance whose value depends on the machine precision.

Denoting $m_o \coloneqq |\Omega|$ and $m_s \coloneqq N (N-1)/2$ 
the number of observation and symmetry constraints, respectively, 
there are $n \coloneqq N^2$ variables and $m \coloneqq m_o + m_s$ constraints 
in \eqref{eq:matrixcompletionminrank}.
We reformulate the model in the form of \eqref{eq:P} by setting 
$f \coloneqq 0$, $D \coloneqq \{0\}$ and a constraint function 
$\func{c}{\R^{N \times N}}{\R^m}$ returning the observation and symmetry constraints 
stacked in vector form.

For our experiments,
we chose $N \in \{ 10, 20 \}$, $\ell = 5$, 
$m_o = \lfloor (n - m_s)/3 \rfloor$, $p = 0.5$ and consider $30$ randomly generated
instances for each value of $N$.
We generate $X \in \R^{N \times \ell}$ by sampling the standard normal distribution, \ie,
$X_{ij} \sim \mathcal{N}(0,1)$, $(i,j)\in\mathcal{I}^2$, and then compute $\Delta$.
Finally, we sample observations by selecting $m_o$ different entries of $\Delta$ 
with uniform probability.

We run our solver \alps{} with default options, 
and abstain from setting an iteration limit for the subproblem solver.
The initial guess $B^0 \in \R^{N \times N}$ is chosen randomly
based on $B_{ij}^0 \sim \mathcal{N}(0,1)$, $(i,j)\in\mathcal{I}^2$,
whereas the dual initial guess is fixed to $y^0 \coloneqq 0$.
We invoke \alps{} directly for solving \eqref{eq:matrixcompletionminrank} 
with the different cost functions mentioned above.
Additionally, the solutions obtained with nuclear norm and Schatten quasi-norm as cost functions,
which are at least continuous,
are used as initial guesses for another round of minimization 
exploiting the discontinuous rank functional.

\begin{figure}[tb]
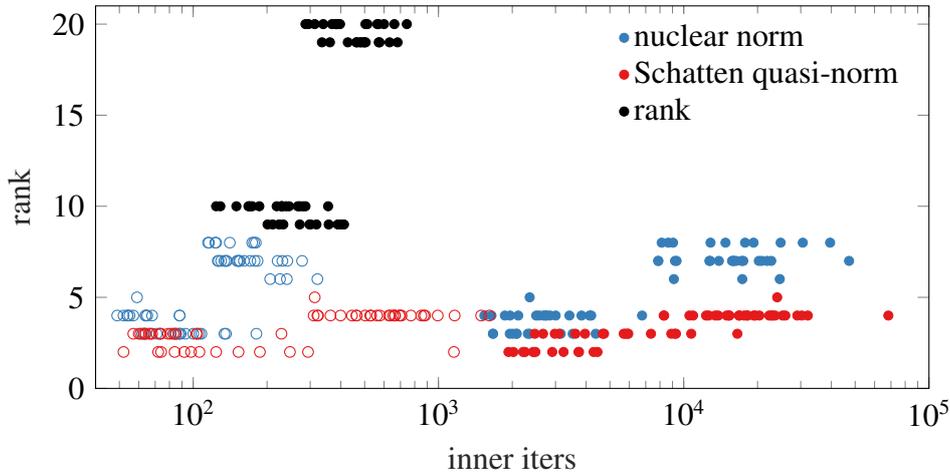

	\centering
	\includetikz{dmc}%
	\caption{%
		Results for the matrix recovery problem \eqref{eq:matrixcompletionminrank}:
		Comparison of (accumulated) inner iteration numbers and rank of the solutions 
		found with different formulations, 
		including warm-started rank minimization (circle).%
	}%
	\label{fig:dmc}%
\end{figure}

We depict the results of our experiments in \cref{fig:dmc}.
Minimization based on the (convex) nuclear norm produces matrices with
rank between $3$ and $8$, while the use of the Schatten quasi-norm 
culminates in solutions having rank between $2$ and $5$.
These findings outperform the direct minimization of the rank which
results in matrices of rank between $9$ and $20$.
This behavior is not surprising since \eqref{eq:matrixcompletionminrank}
possesses plenty of non-global minimizers in case where minimization of the 
discontinuous rank is considered, and \alps{} can terminate in such solutions.
Let us mention that, out of $60$ instances, the warm-started rank minimization 
yields further reduction of the rank in one case after minimization of the Schatten 
quasi-norm and $11$ cases after minimization of the nuclear norm;
in all other cases, no deterioration has been observed.
In summary, \alps{} manages to find feasible solutions 
of \eqref{eq:matrixcompletionminrank} in all cases,
and with adequate objective value in cases where we minimize the
nuclear norm or the Schatten quasi-norm. These solutions can be
used as initial guesses for a warm-started minimization of the rank via
\alps{} or tailored mixed-integer numerical methods.

\section{Conclusions}\label{sec:Conclusions}
	We presented the class of constrained composite optimization problems and proposed a general-purpose solver 
based on an augmented Lagrangian method.
The (outer) augmented Lagrangian loop generates a sequence of subproblems, each one being a dual proximal regularization of the original, 
that can be solved, e.g., by off-the-shelf proximal algorithms for composite optimization.
Requiring only first-order primitives, such as gradient and proximal mapping oracles, and projections onto the constraint set, 
the method is matrix-free and allows the seamless integration of routines for special problem structures.
The proposed method is easily warm started to reduce the number of iterations and can take advantage of accelerated methods.

We have implemented our algorithm in the open-source Augmented Lagrangian Solver (\alps{}), disentangled from modeling tools and subproblem solvers.
Thanks to its low memory footprint and simple, yet fast and robust iterations, \alps{} can handle large-scale problems and is suitable for embedded applications.
We tested our approach numerically with problems arising in mixed-integer optimal control, sparse portfolio optimization and minimum-rank matrix completion.
Illustrative examples showed the flexibility and descriptive power of constrained composite programs 
and the impact of accelerated methods for solving the inner problems.

\subsection*{Acknowledgements}
\TheAcknowledgements

\medskip

\noindent\TheFunding

\clearpage
\phantomsection
\addcontentsline{toc}{section}{References}%
{%
	\small
	\bibliographystyle{plain}
	\bibliography{TeX/biblio}
}

\clearpage
\appendix
\section{On the subproblem solver}\label{sec:appendix}
	In this appendix, we briefly describe the algorithm \panoc{+} from \cite{demarchi2022proximal}, which is
used as a subproblem solver in \cref{alg:ALM}, and discuss some of its properties.

Let us consider the abstract unconstrained, composite optimization problem
\begin{equation}\label{eq:composite_problem}\tag{Q}
	\minimize_{z\in\R^p}
    {}\quad{}
    \omega(z) \coloneqq \varphi(z) + \psi(z)
\end{equation}
under the following standing assumption.
\begin{mybox}
	\begin{ass}\label{ass:PANOC}
		The following hold in \eqref{eq:composite_problem}:
		\begin{enumerate}
			\item\label{ass:phi}%
				$\func{\varphi}{\ZZ}{\R}$ is continuously differentiable with locally Lipschitz continuous gradient;
			\item\label{ass:psi}%
				$\func{\psi}{\ZZ}{\Rinf}$ is proper, lower semicontinuous and prox-bounded with threshold $\gamma_\psi>0$;
			\item\label{ass:inf_omega}%
				$\inf_{z\in\ZZ}\omega(z)>-\infty$.
		\end{enumerate}
	\end{ass}
\end{mybox}

For simplicity of notation, we introduce a set-valued mapping $\ffunc{\FBT_{\gamma}}{\ZZ}{\ZZ}$ for arbitrary $\gamma\in(0,\gamma_\psi)$ by means of
\begin{equation}\label{eq:prox_map}
	\FBT_{\gamma}(z)
	\coloneqq
	\prox_{\gamma\psi}(z-\gamma\nabla\varphi(z)).
\end{equation}
Furthermore, the algorithm makes use of the so-called \emph{forward-backward envelope} (FBE) relative to \eqref{eq:composite_problem}
with stepsize $\gamma\in(0,\gamma_\psi)$ given by
\[
	\omega^{\textup{FB}}_{\gamma}(z)
	\coloneqq
	\min\limits_{w\in\ZZ}
		\varphi(z)+\innprod{\nabla \varphi(z)}{w-z}+\psi(w)+\tfrac1{2\gamma}\|w-z\|^2.
\]
Clearly, for any $\bar z\in\FBT_\gamma(z)$, we have
\begin{equation}\label{eq:PANOC+:fbe}
	\omega^{\textup{FB}}_{\gamma}(z)
	=
	\varphi(z)+\innprod{\nabla \varphi(z)}{\bar z-z}+\psi(\bar z)+\tfrac1{2\gamma}\|\bar z-z\|^2.
\end{equation}

In \cref{alg:PANOC+}, we provide the pseudo code for \panoc{+},
whose peculiarity is the intricate structure emerging at \cref{state:PANOC+:gammaLS,state:PANOC+:tauLS}.
The two backtracking linesearches are entangled,
concurrently affecting both the direction stepsize $\tau_k$ and the proximal stepsize $\gamma_k$.
These persistent adjustments allow \panoc{+}
to construct a tighter merit function $\omega^{\textup{FB}}_{\gamma}$
that better captures the (local) landscape of $\omega$,
obviating the need for global Lipschitz gradient continuity
for the smooth term in \eqref{eq:composite_problem}.

\begin{algorithm}[tbh]
	\caption{PANOC$^+$ \cite{demarchi2022proximal}}%
	\label{alg:PANOC+}%
	\begin{algorithmic}[1]%
\linespread{1.2}\selectfont%
\Require
	\(z^0 \in \ZZ\),
	\(\gamma_0 \in (0,\gamma_\psi)\),
	\(\Delta\geq 0\),
	\(\alpha,\beta\in(0,1)\),
	\(\varepsilon>0\)%
\Initialize
	\(k\gets0\), and start from \cref{state:PANOC+:barz}
	\State
		\label{state:PANOC+:init}%
		\(\gamma_k\gets\gamma_{k-1}\)
	\State
		\label{state:PANOC+:d}%
		Select an update direction \(d^k\in\ZZ\) with \(\|d^k\|\leq \Delta\|\bar z^{k-1}-z^{k-1}\|\) and set \(\tau_k=1\)%
	\State
		\label{state:PANOC+:z+}%
		Set \(z^k=(1-\tau_k)\bar z^{k-1}+\tau_k(z^{k-1}+d^k)\)
	\State
		\label{state:PANOC+:barz}%
		Compute $\bar{z}^k \in \FBT_{\gamma_k}(z^k)$ and set $\Phi_k:=\omega^{\textup{FB}}_{\gamma_k}(z^k)$ as in \eqref{eq:PANOC+:fbe}%
	\If{$\varphi(\bar{z}^k) > \varphi(z^k) + \innprod{\nabla\varphi(z^k)}{\bar{z}^k - z^k} + \tfrac{\alpha}{2\gamma_k}\|\bar{z}^k - z^k\|^2$}%
	\label{state:PANOC+:gammaLS}%
		\Statex*%
			$\gamma_k \gets \gamma_k/2$, and go back to \cref{state:PANOC+:d} if \(k>0\), or \cref{state:PANOC+:barz} if \(k=0\)
	\EndIf
	\If{$\| \tfrac{1}{\gamma_k} (\bar{z}^k - z^k) - \nabla \varphi(\bar{z}^k) + \nabla \varphi(z^k) \| \leq \varepsilon$ }
		\Statex*
			\algfont{Return} $\bar{z}^k$
	\EndIf
	\If{$k>0$ ~{\sc and}~ $\Phi_k > \Phi_{k-1} - \beta\tfrac{1-\alpha}{2\gamma_{k-1}} \|\bar z^{k-1}-z^{k-1}\|^2$}\label{state:PANOC+:tauLS}%
		\Statex*
			$\tau_k \gets \tau_k/2$ and go back to \cref{state:PANOC+:z+}
	\EndIf
	\State
		\(k\gets k+1\) and start the next iteration at \cref{state:PANOC+:init}
\end{algorithmic}
\end{algorithm}

The analysis in \cite{demarchi2022proximal} provides global convergence guarantees for \panoc{+}
under \cref{ass:PANOC}.
Let us recall the basic result associated with \cref{alg:PANOC+} 
that is important in the context of \cref{alg:ALM}.
For the reader's convenience, we present a brief proof of the result as
it is not explicitly stated in \cite{demarchi2022proximal}.
\begin{mybox}
	\begin{prop}\label{prop:PANOC}
		Let $\{z^k\}$ and $\{\bar z^k\}$ be sequences generated by \cref{alg:PANOC+}.
		Furthermore, let $z^\ast$ be an accumulation point of $\{z^k\}$
		and $\{z^k\}_K$ a subsequence such that $z^k\to_Kz^\ast$.
		Then, $z^\ast$ is a stationary point of $\omega$.
		Additionally, $\bar z^k\to_K z^\ast$ holds,
		and for each $\varepsilon>0$ and any large enough $k\in K$,
		$\bar z^k$ is an $\varepsilon$-stationary point of $\omega$.
	\end{prop}
\end{mybox}
\begin{proof}
	Owing to \cite[Thm~4.3]{demarchi2022proximal}, we have
	$\bar z^k\to_K z^\ast$, and $\gamma_k=\gamma$ holds for some
	$\gamma>0$ and large enough $k\in K$.
	Furthermore, this result gives boundedness of the expressions
	\begin{align*}
		\Phi_k
		\coloneqq
		\varphi(z^k)+\innprod{\nabla\varphi(z^k)}{\bar z^k-z^k}+\psi(\bar z^k)
		+
		\tfrac{1}{2\gamma_k}\|\bar z^k-z^k\|^2,
	\end{align*}
	so that taking the lower limit $k\to_K\infty$ yields $z^\ast\in\dom\psi$.
	Next, \cref{state:PANOC+:barz} of \cref{alg:PANOC+} yields
	\begin{align*}
		\omega(z^\ast)
		&\leq
		\liminf\limits_{k\to_K\infty}\Phi_k
		\\
		&\leq
		\liminf\limits_{k\to_K\infty}
		\left(
			\varphi(z^k)+\innprod{\nabla\varphi(z^k)}{z^\ast-z^k}+\psi(z^\ast)
			+
			\tfrac{1}{2\gamma_k}\|z^\ast-z^k\|^2
		\right)
		\\
		&\leq
		\limsup\limits_{k\to_K\infty}
		\left(
			\varphi(z^k)+\innprod{\nabla\varphi(z^k)}{z^\ast-z^k}+\psi(z^\ast)
			+
			\tfrac{1}{2\gamma_k}\|z^\ast-z^k\|^2
		\right)
		\\
		&=
		\omega(z^\ast),
	\end{align*}
	giving $\bar z^k\toattentive{\omega}_K z^\ast$ by continuity of $\varphi$.
	Considering the stationarity condition resulting from evaluation of the
	proximal map $\FBT_{\gamma_k}$, 
	\[
		0\in\nabla \varphi(z^k)+\partial\psi(\bar z^k)+\tfrac1{\gamma_k}(\bar z^k-z^k)
	\]
	holds for each $k\in K$, giving
	\[
		\tfrac{1}{\gamma_k}(z^k-\bar z^k)
		+
		\nabla\varphi(\bar z^k)-\nabla\varphi(z^k)
		\in 
		\nabla\varphi(\bar z^k)+\partial \psi(\bar z^k)
		=
		\partial\omega(\bar z^k).
	\]
	Taking the limit $k\to_K\infty$ while respecting continuous differentiablity of
	$\varphi$, the result follows.
\end{proof}

Let us mention that slightly weaker convergence guarantees can be obtained
for \panoc{+} whenever the evaluation of the proximal mapping $\FBT_{\gamma_k}$ 
in \cref{state:PANOC+:barz} of \cref{alg:PANOC+}
is done inexactly, see \cite[\S 4]{demarchi2022proximal} for details.

Finally, in light of \cref{sec:numresults}, we shall comment on the acceleration mechanism in \panoc{+}.
Although robust to arbitrary choices of (bounded) 
directions $d^k$, the practical performance of \cref{alg:PANOC+} is strongly affected by
the specific selection; we refer to \cite[\S 4.3]{themelis2018proximal} 
for an overview on some potential update directions.
In the numerical experiments, we consider two strategies
for executing \cref{state:PANOC+:d} of \cref{alg:PANOC+}.
First, we may select $d^k \coloneqq \bar{z}^{k-1} - z^{k-1}$, 
so that $z^k = \bar{z}^{k-1}$, effectively reducing the algorithm to
an adaptive proximal gradient method,
without any acceleration \cite[\S 4.4]{demarchi2022proximal}.
Second, as a baseline, we use the default acceleration strategy in \ProximalAlgorithmsLink{}, 
namely LBFGS directions with memory $5$.
Inspired by quasi-Newton methods, 
these are recursively constructed by keeping memory 
of pairs $z^{k+1} - z^k$ 
and $r^{k+1} - r^k$, with $r^k \coloneqq z^k - \bar{z}^k$, 
and retrieving $d^k \coloneqq - H^k r^k$ 
by simply performing scalar products \cite{liu1989limited}.
Herein, the linear operator $H_k$ mimics the (inverse) fixed-point residual mapping 
associated to the splitting scheme in a neighborhood of $z^k$
\cite{stella2017simple,themelis2018forward}.
Notice that, as the geometry of the residual mapping depends on the proximal stepsize, 
(the memory of) the LBFGS approximation is reset every time 
the stepsize is adapted \cite[\S 3.1]{demarchi2022proximal}.

\end{document}